\documentclass[12pt]{amsart}
\usepackage{amsmath,amssymb,amsbsy,amsfonts,amsthm,latexsym, amsopn,amstext,amsxtra,euscript,amscd,mathrsfs,color,bm, cite}
\usepackage{todonotes}   
\usepackage{url}
\usepackage[colorlinks,linkcolor=blue,anchorcolor=blue,citecolor=blue]{hyperref}
\usepackage{color}
\usepackage{comment}

\newtheorem{theorem}{Theorem}[section]
\newtheorem{lemma}[theorem]{Lemma}
\newtheorem{crit}[theorem]{Criterion}

\newtheorem{cor}[theorem]{Corollary}
\newtheorem{prop}[theorem]{Proposition}

\newtheorem{rem}[theorem]{Remark}
\newtheorem*{Heilbronn's Criterion}{Heilbronn's Criterion}

\begin{document}

\def\cA{{\mathcal A}}
\def\cB{{\mathcal B}}
\def\cC{{\mathcal C}}
\def\cD{{\mathcal D}}
\def\cE{{\mathcal E}}
\def\cF{{\mathcal F}}
\def\cG{{\mathcal G}}
\def\cH{{\mathcal H}}
\def\cI{{\mathcal I}}
\def\cJ{{\mathcal J}}
\def\cK{{\mathcal K}}
\def\cL{{\mathcal L}}
\def\cM{{\mathcal M}}
\def\cN{{\mathcal N}}
\def\cO{{\mathcal O}}
\def\cP{{\mathcal P}}
\def\cQ{{\mathcal Q}}
\def\cR{{\mathcal R}}
\def\cS{{\mathcal S}}
\def\cT{{\mathcal T}}
\def\cU{{\mathcal U}}
\def\cV{{\mathcal V}}
\def\cW{{\mathcal W}}
\def\cX{{\mathcal X}}
\def\cY{{\mathcal Y}}
\def\cZ{{\mathcal Z}}

\def\A{{\mathbb A}}
\def\B{{\mathbb B}}
\def\C{{\mathbb C}}
\def\D{{\mathbb D}}
\def\E{{\mathbb E}}
\def\F{{\mathbb F}}
\def\G{{\mathbb G}}
\def\I{{\mathbb I}}
\def\J{{\mathbb J}}
\def\K{{\mathbb K}}
\def\L{{\mathbb L}}
\def\M{{\mathbb M}}
\def\N{{\mathbb N}}
\def\P{{\mathbb P}}
\def\Q{{\mathbb Q}}
\def\R{{\mathbb R}}
\def\S{{\mathbb S}}
\def\T{{\mathbb T}}
\def\U{{\mathbb U}}
\def\V{{\mathbb V}}
\def\W{{\mathbb W}}
\def\X{{\mathbb X}}
\def\Y{{\mathbb Y}}
\def\Z{{\mathbb Z}}

\def\fl#1{\left\lfloor#1\right\rfloor}
\def\rf#1{\left\lceil#1\right\rceil}
\def\eps{\varepsilon}
\def\mand{\qquad\mbox{and}\qquad}

\def\sssum{\mathop{\sum\ \sum\ \sum}}
\def\ssum{\mathop{\sum\, \sum}}
\def\ssumw{\mathop{\sum\qquad \sum}}

\def\inv#1{\overline{#1}}
\def\num#1{\mathrm{num}(#1)}
\def\dist{\mathrm{dist}}

\def\fA{{\mathfrak A}}
\def\fB{{\mathfrak B}}
\def\fC{{\mathfrak C}}
\def\fU{{\mathfrak U}}
\def\fV{{\mathfrak V}}

\newcommand{\bflambda}{{\boldsymbol{\lambda}}}
\newcommand{\bfxi}{{\boldsymbol{\xi}}}
\newcommand{\bfrho}{{\boldsymbol{\rho}}}
\newcommand{\bfnu}{{\boldsymbol{\nu}}}

\def\GL{\mathrm{GL}}
\def\SL{\mathrm{SL}}

\def\Hba{\overline{\cH}_{a,m}}
\def\Hta{\widetilde{\cH}_{a,m}}
\def\Hb1{\overline{\cH}_{m}}
\def\Ht1{\widetilde{\cH}_{m}}

\def\flp#1{{\left\langle#1\right\rangle}_p}
\def\flm#1{{\left\langle#1\right\rangle}_m}

\def\Zm{\Z/m\Z}

\def\Err{{\mathbf{E}}}
\def\O{\mathcal{O}}

\def\cc#1{\textcolor{red}{#1}}

\newcommand{\km}[1]{{\color{purple}{KM: #1}} }
\newcommand{\ab}[1]{{\color{blue}{AB: #1}} }
\newcommand{\bk}[1]{{\color{red}{BK: #1}} }
\newcommand{\bg}[1]{{\color{magenta}{BG: #1}} }
\newcommand{\vs}[1]{{\color{teal}{VS: #1}} }

\newcommand{\comm}[1]{\marginpar{%
\vskip-\baselineskip 
\raggedright\footnotesize
\itshape\hrule\smallskip#1\par\smallskip\hrule}}

\def\xxx{\vskip5pt\hrule\vskip5pt}

\def\dmod#1{\,\left(\textnormal{mod }{#1}\right)}

\allowdisplaybreaks


\title{The determination of norm-Euclidean cyclic cubic fields}

\author[G. Bagger]{Gustav Kj\ae rbye Bagger}\thanks{}
\address{School of Science, The University of New South Wales Canberra, Australia}
\email{g.bagger@unsw.edu.au}

\author[A. Booker]{Andrew R. Booker}\thanks{}
\address{School of Mathematics, University of Bristol, England}
\email{andrew.booker@bristol.ac.uk}

\author[B. Kerr]{Bryce Kerr}\thanks{BK is partially supported by ARC Grants DE220100859 and DP230100534.}
\address{School of Science, The University of New South Wales Canberra, Australia}
\email{bryce.kerr@unsw.edu.au}

\author[K. McGown]{Kevin J.~McGown}\thanks{}
\address{Department of Mathematics and Statistics\\ California State University, Chico, USA}
\email{kmcgown@csuchico.edu }

\author[V. Starichkova]{Valeriia Starichkova}\thanks{VS is partially supported by ARC Grant DP240100186 and the Australian Mathematical Society Lift-off Fellowship.}
\address{School of Science, The University of New South Wales Canberra, Australia}
\email{v.starichkova@unsw.edu.au}

\author[T. Trudgian]{Tim Trudgian}\thanks{TT is partially supported by ARC Grants FT160100094 and  DP240100186.}
\address{School of Science, The University of New South Wales Canberra, Australia}
\email{timothy.trudgian@unsw.edu.au}

\date{\today}
\pagenumbering{arabic}


\begin{abstract}
It is known on the Generalised Riemann Hypothesis that there are precisely $13$ cyclic cubic fields that are norm-Euclidean. Unconditionally, there is a gap between analytic estimates which hold for all sufficiently large conductors and computational techniques. In this paper, we establish new results concerning explicit bounds for cubic non-residues and refine previous computational techniques, enabling us to completely characterise all norm-Euclidean cyclic cubic fields.
\end{abstract}

\maketitle

\section{Introduction}\label{Intro}
Let $K$ be a number field.  Let $\O_K$ denote the ring of integers and $N$ denote the norm map.
We say that $K$ is norm-Euclidean if for all $\alpha\in K$ there exists $\beta\in\O_K$
such that $|N(\alpha-\beta)|<1$.
This is the classical generalisation of the Euclidean algorithm to the number field setting,
considered by Gauss, Dirichlet, and Dedekind~\cite{Gauss1, Gauss2, Dirichlet, Dedekind}.
In the more modern vocabulary, this condition is equivalent to saying that $\O_K$ is a Euclidean domain
with respect to the absolute value of the norm.
Number fields which are norm-Euclidean must have class number one, but the converse does not 
hold in general.

When $K$ is quadratic it is known (see \cite{Chat} and~\cite{Barnes}) that
$K=\mathbb{Q}(\sqrt{d})$ with $d$ square-free is norm-Euclidean precisely when
\begin{equation*}
d= -1, \pm 2, \pm 3, \pm 7, \pm 11, 5, 6, 13, 17, 19, 21, 29, 33, 37, 41, 57, 73.
\end{equation*}
The imaginary quadratic case of this determination
was known to Dirichlet and Dedekind. The real quadratic case
is more difficult as the theorem was not proved until around 1950.
Although this is often attributed to Chatland and Davenport~\cite{Chat},
all told it was the culmination of the work of over a
dozen mathematicians, most of which took place
from 1930 to 1950.
The work of Barnes and Swinnerton-Dyer~\cite{Barnes}
brought to light an error in the 
initial classification; it was believed for a time that
$\Q(\sqrt{97})$ was norm-Euclidean.
These final two results occurred over a century after Wantzel~\cite{Want} first showed that
the fields
$\Q(\sqrt{d})$, when $d=2,3,5$, are norm-Euclidean!

As the quadratic case is completely settled, we consider the case of cubic fields.
Davenport \cite{Dav1950} proved that there are only finitely many complex cubic fields that are norm-Euclidean --- see also Lemmermeyer \cite[Sec.\ 5.1]{Lemon}.
This leaves open the case of totally real cubic fields.
Heilbronn \cite{Hei1950} suggested there may be infinitely many norm-Euclidean fields in this situation
by saying that he would be ``surprised to learn that the analogue [of the finiteness theorem] is true in this case''.
In fact, this is still an open problem.

In order to make progress, one may restrict to the case of Galois fields.
Indeed, a very natural question is as follows:
Which Galois number fields of degree $n$ are norm-Euclidean?
Until this point, this question has been answered for $n=2$, but no other value of $n$.
Of course, when $n$ is a prime, there is only one possible Galois group,
the cyclic group of order $n$; hence we may simply speak of cyclic fields.

We specialise to the case of cyclic cubic fields, which are necessarily totally real.
In this case, Heilbronn proved that only finitely many are norm-Euclidean.
By the conductor--discriminant formula, the discriminant of such a field takes the form $\Delta=f^2$,
and by genus theory, it must be that either $f$ is a prime with $f\equiv 1\pmod{3}$ or $f=9$.
The latter fact follows from assuming, as we may, that $K$ has class number one.

Godwin and Smith \cite{Godwin, Smith, GodSmith} showed that the only cyclic cubic fields that are norm-Euclidean
with $f\leq 10^4$ are the following:
\begin{equation}\label{dog}
f = 7,9,13,19,31,37,43,61,67,103,109,127,157.
\end{equation}
At the time, no one seemed to conjecture that this list was complete;
the difficulty was that there was no known upper bound on the discriminant for such a field.
In~\cite{McGown1}, McGown gave the first such upper bound.
He proved (see~\cite{McGown1, McGown2}),
under the Generalised Riemann Hypothesis (GRH),
that the list in (\ref{dog}) is complete, and unconditionally,
that any exception must satisfy $f\in(10^{10}, 10^{70})$.
This range was reduced to $(2\times  10^{14}, 10^{50})$ by Lezowski and McGown \cite{Lezowski}.
The upper bound was reduced slightly by Francis \cite{Forrest} to $3.8\times  10^{49}$.
Our paper completely solves the problem. 
\begin{theorem}\label{T:big}
A cyclic cubic field is norm-Euclidean if and only if it has conductor
\begin{equation*}
f \in \{7,9,13,19,31,37,43,61,67,103,109,127,157\}.
\end{equation*}
\end{theorem}


In the next section, we give a brief sketch of our argument which
serves as an overview of the paper.

\section{Overview}

In this section, we outline the strategy for
showing that no norm-Euclidean cyclic cubic field has conductor 
$f\ge 2\times 10^{14}$.

Let $K$ be a cyclic cubic field.  As discussed in Section~\ref{Intro},
we may assume that $K$ has class number one.  In particular,
$\Delta_K=f^2$ where $f$ is a prime with $f\equiv 1\pmod{3}$,
after disposing of the single field with $f=9$.
We will denote by $\mathcal{N}$
the image of $\O_K$ under the norm $N$; i.e.\, $\mathcal{N} = N(\O_K)$.
The following result is key for our argument:

\begin{Heilbronn's Criterion}
If $f$ can be represented as $f = a + b$,
where $a,b>0$,
$a,b\not\in\mathcal{N}$, and
$a$ is a cubic residue modulo $f$,
then $K$ is not norm-Euclidean.
\end{Heilbronn's Criterion}

This result was established by Heilbronn~\cite{Hei1950} who proved the finitude of cyclic cubic norm-Euclidean number fields $K$ by showing that the criterion above holds when $f$ is sufficiently large.

Let $q_1<q_2$ be the two smallest inert primes in $K$.  This is equivalent to saying $q_1$ and $q_2$ are the smallest prime cubic non-residues modulo $f$, meaning that for a primitive character $\chi$ modulo $f$ of order $3$ (defined uniquely up to the complex conjugation), $\chi(q_1), \chi(q_2) \not \in \{0,1\}$.
The idea to find suitable $a,b$ is to seek an expression of the form $f=uq_1+vq_2$ for suitable $u,v$.
Indeed,
given $n\in \Z$, one has $n\not\in\mathcal{N}$ if and only if $3\nmid\nu_p(n)$
for some inert prime $p$;
here $\nu_p(n)$ denotes the $p$-adic valuation of $n$, i.e.,
the highest power of $p$ dividing $n$.

Let $m$ be the smallest positive integer coprime to $q_1 q_2$ such that $\chi(m) = \chi^{-1}(q_2)$. It is explained in the proof of\cite[Proposition 6.1]{Lezowski} why
suitable $a,b$ exist provided, say, $q_1\neq 2$ and $3 q_1 q_2 m\leq f$.
Essentially,
Heilbronn's criterion will hold for a given $f$ provided we can show an estimate of the form 
\begin{align}
\label{eq:heil}
q_1 q_2 m \ll f
\,,
\end{align}
with a small enough implicit constant.
We refer the reader to Criterion \ref{criterion: expl-non-Eucl}
for a more nuanced version of (\ref{eq:heil}) with all the details rendered visible.
In this way, the question of non-Euclideanity reduces to giving upper bounds for
$q_1$, $q_2$, and $m$.

Since cubic non-residues may be detected via multiplicative characters, establishing~\eqref{eq:heil}
is closely related to obtaining
bounds on character sums.
The most natural such character sum is
\begin{align}
\label{eq:charsum-overview}
\sum_{\substack{n\le x}}\chi(n) 
\,,
\end{align}
but we will also require the use of 
more complicated but related character sums
such as
\begin{align}
\label{eq:charsum-overview2}
  \sum_{x=1}^f\left|\sum_{k=0}^{h-1}\chi(x+k)\right|^{2r}
  \,,\qquad
  \sum_{1\leq t\leq T}
  \sum_{1\leq u\leq U}
  \sum_{1\leq v\leq V}
  \chi(v-ut)
  \,.
\end{align}

Under GRH, we may estimate~\eqref{eq:charsum-overview} in the range $x\gg\log^2{f}$, leading to the inexplicit bound $q_1 = O(\log^2 f)$ which was first established by Ankeny \cite{Ankeny} (and later generalised by Lagarias--Montgomery--Odlyzko \cite{Lagarias}) with the corresponding explicit bound given by Bach~\cite{Bach}. The argument from~\cite{Bach} was generalised in \cite{McGown2} by McGown to obtain an explicit bound of the same magnitude for $q_1$, $q_2$, and $m$, and thus~\eqref{eq:heil} holds which proves Theorem \ref{T:big} under GRH.

Establishing~\eqref{eq:heil} in an explicit form
without recourse to additional hypotheses is difficult.
The best known character sum estimates for~\eqref{eq:charsum-overview} are due to Burgess~ \cite{Burgess} and hold in the range
$x\gg_\varepsilon f^{1/4+\varepsilon}$.
Using this, one can show $q_1, q_2, m\ll_\varepsilon f^{1/4+\varepsilon}$,
and hence
\eqref{eq:heil} holds when $f$ is sufficiently large.
Combining Burgess' estimate with techniques from sieve theory (often referred to as Vinogradov's trick in the context of character sums),
it is possible to further show 
that $q_1, q_2\ll_\varepsilon f^{1/(4\sqrt{e}) + \varepsilon}$
(see for example~\cite{Pollack}).

Using state-of-the-art bounds for non-residues (see~\cite{Norton, Trev2, Ma, McGown3, Lezowski})
and explicit versions of Burgess' character sum estimate
(see~\cite{McGown1,Trev,Forrest}),
one can show that Heilbronn's criterion holds when 
$f\geq 10^{50}$.  The goal of this paper is to cover the remaining values of $f$.
The main novelty of this paper is to provide refinements of various character sum
and sieve arguments using a variety of methods, and apply them on a case-by-case
basis depending on the relative sizes of $q_1,q_2$ to cover the remaining range.
The range of $10^{22}\leq f\leq 10^{50}$ is covered by our analytic estimates
(see Sections~\ref{sec:I}-\ref{sec:II})
and the range $2\times 10^{14}\leq f\leq 10^{22}$ is covered by our 
computational method (see Section~\ref{numerics}).



Our first bound on cubic non-residues
is Theorem~\ref{T:nonresidues}, which is proved in Section~\ref{sec:I}.  This
is based on techniques of McGown--Trudgian~\cite{KevTim},
Ma--McGown--Rhodes--Wanner~\cite{Ma},
and
takes advantage of the cubic nature of the character.
As is true for many character sum estimates that ultimately invoke the Weil bound,
the proof of this result involves considering the first character sum 
on the left-hand side of~\eqref{eq:charsum-overview2}; for this,
we make use of a formula of Trevi\~no (see~\cite{Trev}).
As a consequence of Theorem~\ref{T:nonresidues}, supplemented by a P\'olya--Vinogradov bound (Lemma~\ref{lem:pv}) when $q_1q_2$ is small, we deduce Theorem~\ref{T:1} which states
that Heilbronn's condition holds whenever $q_1\in\{2,3,5,7\}$ and $f\geq 10^{14}$.
In a similar manner, for any fixed $q_1$
we obtain such a bound on~$f$;
these are tabulated in
Table~\ref{tab:q1thresh} for $q_1\leq 199$.  For example, this result says that Heilbronn's
condition holds whenever $q_1\leq 199$ and $f\geq 1.60\times 10^{19}$.

An important reduction is provided by
Proposition~\ref{fixedprop14}, which occurs in Section~\ref{sec:mcg}.
The idea here is that character sum estimates of McGown~\cite{McGown1}
(see Lemma~\ref{lem:mcgownR}),
supplemented with our new
Theorem~\ref{T:nonresidues},
allow us to cover smaller sized values of $q_1$ or $q_2$.
When $f\geq 10^{22}$,
all this allows us to reduce to the situation where
\begin{align}
\label{eq:range1}
q_1\ge  379 \quad \text{and} \quad q_2\ge 1.3f^{1/6}.
\end{align}

Dealing with this range of $q_1,q_2$ is the focus of Section~\ref{charsum} and Section~\ref{sec:vino}. In Section~\ref{charsum} we prove some new explicit character sum estimates. Our main result of Section~\ref{charsum} is Theorem~\ref{thm:1}, which significantly improves current explicit estimates for
the second character sum in~\eqref{eq:charsum-overview2}.
One of the key steps in the traditional Burgess--Karatsuba amplification is a bound on the multiplicative energy of an interval:
\begin{align}
\label{eq:gcd}
\#\{ 1\le n_i\le N, \ \ 1\le u_i \le U \ : \ n_1u_1\equiv n_2u_2 \pmod{f} \}\approx NU(\log{U}).
\end{align}
The $\log{U}$ comes about from $\gcd$ sum estimates and has an unacceptable contribution if
we are aiming to cover the range $f \ge 10^{22}$. Since~\eqref{eq:gcd}
can be expressed as the average lattice point counting problem
\begin{align*}
\sum_{\lambda \hspace{-1ex}\pmod{f}}I(\lambda)^2
\end{align*}
where
\begin{align}
\label{eq:il}
I(\lambda)=\#\{ 1\le n \le N, \ \ 1\le u \le U \ : \ n\equiv\lambda u \pmod{f} \},
\end{align}
by analysing the smallest vector $(n,u)$ in~\eqref{eq:il} it is possible to show that not many values of $\lambda$ contribute to the $\log{U}$ term in~\eqref{eq:gcd}. Thus we may remove these bad values at an earlier stage of the argument, leading to good numerics in Theorem~\ref{thm:1}.

Our next main result is Corollary~\ref{lem:firsttwo1}
to Lemma~\ref{lem:firsttwo11}
in Section~\ref{sec:vino} which combines our character sum estimates from Section~\ref{charsum} with Vinogradov's trick (see, for example~\cite{Vin}). This can be considered a Lipschitz type principle for partial sums of multiplicative functions where one may relax the upper bound on~\eqref{eq:charsum-overview} with the advantage of extending the range of summation $x$ for which the bound is non-trivial. This requires sieving out $q_1,q_2$ from the character sums occurring in Theorem~\ref{thm:1}.
Our previous reduction~\eqref{eq:range1} is necessary in order to control the numerical constants coming from the sieve.
These techniques combine to establish Theorem~\ref{thm:main1} which shows that if $q_1\ge 379$
and $f\in(10^{22}, 10^{50})$
then
\begin{align}
\label{eq:range2}
 q_1\le q_2\le 36.88f^{0.21769}.
\end{align}
Thus, to establish~\eqref{eq:heil} for
$f$ in this range,
and thereby prove Theorem~\ref{T:2},
it remains to combine~\eqref{eq:range2} with suitable bounds for $m$ which are provided by Proposition~\ref{thm:main2} and Lemma~\ref{lem:main2m}. 

Finally, 
the range $2\times 10^{14}\le f\le 10^{22}$ is verified by the
computational method described in Section~\ref{numerics}.
Our algorithm for computing cubic non-residues is similar to the method of enumerating pseudosquares,
with the main difference being that our sieve region is an ellipse rather than an interval.
Essentially the method computes $q_1$ and applies Theorem~\ref{T:1}, if applicable.
Otherwise, it continues to compute $q_1,q_2,m$
and verifies
Criterion \ref{criterion: expl-non-Eucl} directly.
The outcome of the computation is stated as Theorem~\ref{T:compute}.
See Section~\ref{numerics}
for the full details.

\begin{figure}
\centering
\includegraphics[width=0.85\textwidth]{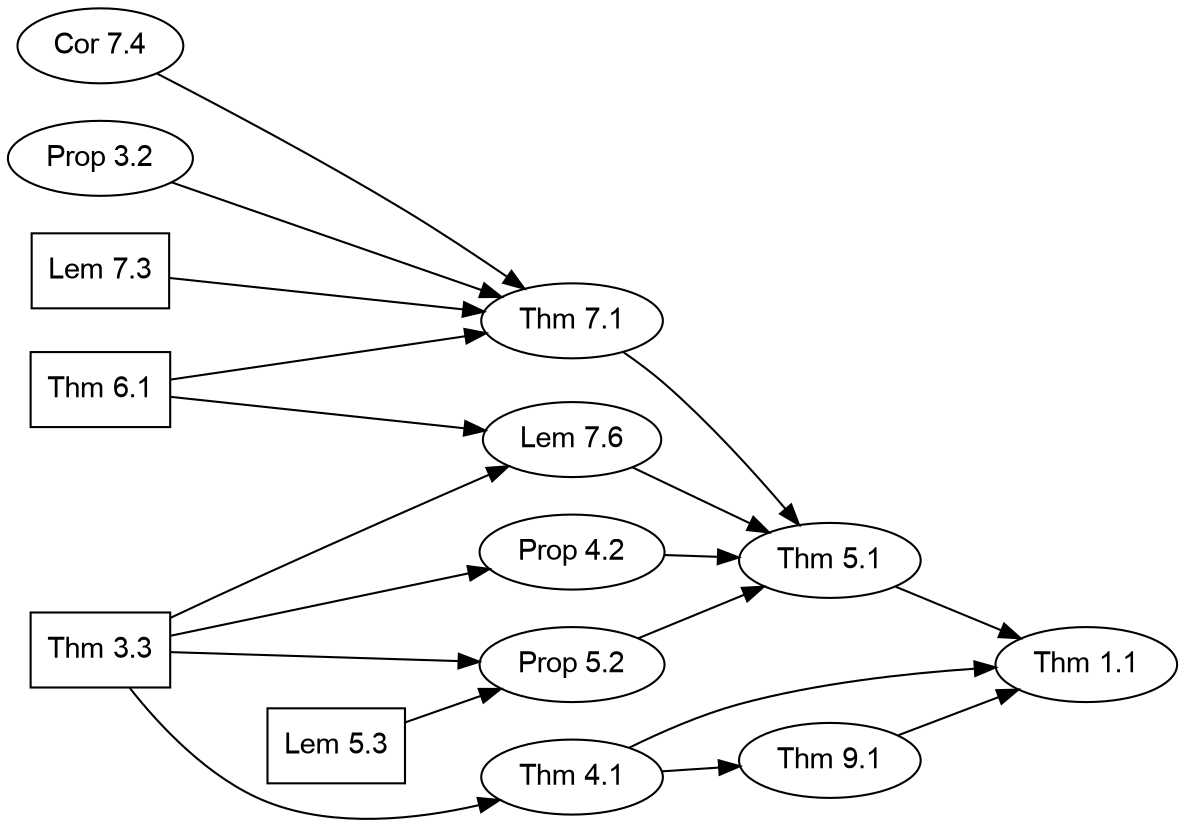}
\caption{Dependencies among the main results.}
\label{fig:dependencies}
\end{figure}

Summarizing, it is fair to say that Theorem~\ref{T:nonresidues}, Theorem~\ref{thm:1}, and Lemma~\ref{lem:firsttwo11}
are our main technical results which, through a series of intermediary results
(including computation), allow us to establish Theorems~\ref{T:1},~\ref{T:2}, and~\ref{T:compute}.
These latter three results taken together prove Theorem~\ref{T:big}.
The logical structure of the proof is illustrated in Figure~\ref{fig:dependencies}.

\section{Cubic non-residues I}
\label{sec:I}

In this section, we give a refinement of the criterion for
the failure of norm-Euclideanity given in~\cite{Lezowski}
(see Criterion~\ref{criterion: expl-non-Eucl})
and establish an explicit bound for cubic non-residues (see Theorem~\ref{T:nonresidues}).

As before, let $K$ be a cyclic cubic field with a prime conductor $f$ and let us fix a Dirichlet character $\chi$ modulo $f$ of order $3$ (there are only two such characters, and they are complex conjugates). An integer $n$ is called a non-residue modulo $f$ if $\chi(n) \not \in \{0,1\}$. As described in~\cite{McGown1}, the norm-Euclideanity of $K$ is connected with the distribution of the first non-residues modulo $f$. Let $q_1<q_2$ denote the smallest prime cubic non-residues modulo $f$.

\begin{crit}\label{criterion: expl-non-Eucl}
\textnormal{Assume there exists some $m \in \N$ satisfying all of
\begin{align*}
&(m,q_1q_2)=1, \quad \chi(m) = \chi^{-1}(q_2), \\
&\text{if } q_1=2: 3q_2 m \leq f, \\
&\text{if } q_1\neq 2: 
\begin{cases}
    \max{\{3q_1 q_2 m,10 q_1^2 q_2\}} \leq f & \text{if} \quad q_2^2 m \equiv f \dmod{q_1} \text{ and } q_2 < 2 q_1, \\
    \max{\{2q_1 q_2 m,10 q_1^2 q_2\}} \leq f & \text{otherwise.}
\end{cases} 
\end{align*}
Then $K$ is not norm-Euclidean.}
\end{crit}
\begin{proof}
The case $q_1 = 2$ follows from Theorem 3.1 and Corollary 5.2 in \cite{McGown1}. The case $q_1 \neq 2$ is a slightly refined version of \cite[Prop.\ 6.1]{Lezowski}: to get this refinement, we note that the condition $3q_1 q_2 m \leq f$ may be replaced by $2q_1 q_2 m \leq f$ in the whole proof of \cite[Prop.\ 6.1]{Lezowski} except the case (d)(i) in which $q_2^2 m \equiv f \dmod{q_1}$ and $q_2 < 2 q_1$.
\end{proof}

This criterion is based on Heilbronn's original criterion \cite{Hei1950} and is the means by which we are able to rule out cyclic cubic fields lacking norm-Euclideanity. In order to invoke Criterion~\ref{criterion: expl-non-Eucl} for a conductor $f$, one constructs a sufficiently small non-residue $m$, which is coprime to $q_1$ and $q_2$. Theorem \ref{T:nonresidues} below ensures the existence of such an $m$ and will be one of our main detection tools for cyclic cubic fields lacking norm-Euclideanity. This theorem relies on upper bounds for the sum
\begin{equation*}
S_\chi(f,h,r):=\sum_{x=1}^{f}\left|\sum_{k=0}^{h-1}\chi(x+k)\right|^{2r}, \quad r,h \in \N
\,.
\end{equation*}

\begin{prop}[\cite{Trev2}, Remark~2.1]\label{prop2}
We have
\begin{equation*}
    S_\chi(f,h,r)\leq f^{1/2}h^{2r}  W(f,h,r),
\end{equation*}
where
\begin{equation*}
  W(f,h,r):=2r-1 + \frac{f^{1/2}}{h^{r}} r! \, d_r(h),
\end{equation*}
and
\begin{equation*}
  d_r(h) = \frac{1}{r! \, h^r} \sum_{n = 0}^{\lfloor \frac{r}{3} \rfloor} \left( \frac{r!}{n!(3!)^n} \right)^2 \frac{h^{r-n}}{(r - 3n)!}.
\end{equation*}
The values of $d_r(h)$, $1 \leq r \leq 6$ are provided in the table below.
\begin{center}
  \begin{tabular}{|c|c|}
  \hline
  $r$ & $d_r(h)$\\
  \hline
  $1$ & $1$\\
  $2$ & $1$\\
  $3$ & $1+1/(6h)$\\
  $4$ & $1+2/(3h)$\\
  $5$ & $1+5/(3h)$\\
  $6$ & $1+10/(3h)+5/(36h^2)$ \\
  \hline
  \end{tabular}
  \end{center}
\end{prop}

\begin{theorem}\label{T:nonresidues}
Let $\chi$ be a cubic Dirichlet character modulo a prime $f$.
Let $\omega$ denote a non-trivial cube root of unity.
Let $r,h,u,v\in\N$ where $u=u_{1}u_{2}\dots u_{k}$ is the product of pairwise distinct primes $u_{i}<h$, and $v=v_{1}v_{2}\dots v_{\ell}$ is the product of pairwise distinct primes satisfying $h\leq v_{i}<f$.
Let $H \in (0,f)$ be a real number and set $X=H/h$.  
Suppose $r\leq 9h$, $3\ell<h$, $\frac{X}{u}\geq 1$, $X\geq 2$, and $2HX\leq f$, and let
\begin{equation}
E_u(X) = 1 - \frac{\pi^2}{6}\sigma(u) \left( \frac{\sigma(u)}{4} + \frac{\phi(u)}{u} + \frac{\phi(u)}{X} \right)\frac{1}{X},
\end{equation}
where $\sigma(u) = \sum_{d \mid u} d$; for $u=1$, the term $\frac{\phi(u)}{u}$ here is to be replaced by $2$ (see the remark following Lemma~\ref{L:sumstuff}). If $E_u(X) > 0$ and
\begin{equation} \label{eq: THm6-main-condition}
    \frac{1}{E_u(X)}\frac{\pi^2}{6}\frac{\sigma(u)}{\phi(u)}uhf^{1/2}\left(\frac{2h}{h-3 \ell}\right)^{2r}\frac{W(f,h,r)}{H^2}<1,
\end{equation}
then there exists a positive integer $m \leq H$ satisfying $(m,uv)=1$ and $\chi(m)=\omega$.
\end{theorem}
In most of our applications, $uv = q_1q_2$ and $h \in \N$ will determine whether $q_1,q_2$ are factors of $u$ or $v$. The value of $H$ will be chosen so that $m\le H$ aligns with Criterion~\ref{criterion: expl-non-Eucl}. 

We next state some preliminary results required for the proof of Theorem \ref{T:nonresidues}. For $t,q \in \Z$ with $0\leq t<q\leq X$ and $(t,q)=1$, we define the intervals
$$ \mathcal{I}(q,t):=\left(\frac{tf}{q},\,\frac{tf+H}{q}-h+1\right]\, \quad \text{and} \quad
\mathcal{J}(q,t):=\left[\frac{tf-H}{q},\,\frac{tf}{q}-h+1\right)\,. $$

\begin{lemma} \label{IJintervals}
Let us keep the notations from Theorem \ref{T:nonresidues} and suppose $0 \leq t < q \leq X$ with $(t,q) = 1$. 
\begin{enumerate}
\item 
Let $0\leq n\leq h-1$.
If $z\in\mathcal{I}(q,t)$ then $q(z+n)-ft\in(0,H]$,
and if $z\in\mathcal{J}(q,t)$ then $q(z+n)-ft\in[-H,0)$.
\item
If $X\geq 2$ and $2HX \leq f$, then $\mathcal{I}(q,t)$, $\mathcal{J}(q,t)$, are pairwise disjoint subsets of $[-H,f-H)$.
\end{enumerate}
\end{lemma}

\begin{proof}
See \cite[Proposition 3]{KevTim}.
\end{proof}

For $u \in \N$, let $N_u(X)$ denote the number of integers in the union of all the intervals $\mathcal{I}(q,t)$ and $\mathcal{J}(q,t)$ satisfying $u\parallel q$, namely
\begin{equation} \label{eq: def-N-u}
    N_u(X):=\sum_{\substack{0\leq t<q\leq X\\(t,q)=1\\u\parallel q}}\sum_{z\in\mathcal{I}(q,t)\sqcup\mathcal{J}(q,t)}1,
\end{equation}
where $u\parallel q$ should be interpreted as $u\mid q$ and $(q/u,u)=1$.

In the next two lemmas, $\vartheta$ will denote a real number with $\vert \vartheta \vert \leq 1$; $\vartheta$ does not need to be the same at each appearance.
\begin{lemma}\label{bowmore}
For a real number $X\geq 1$ and an integer $u>1$,
\begin{equation}\label{lagavulin}
X\sum_{\substack{n\leq X\\(n,u)=1}}1-\sum_{\substack{n\leq X\\(n,u)=1}}n
=\frac{1}{2}\frac{\phi(u)}{u}X^2+\frac{\vartheta}{8}\sigma(u).
\end{equation}
\end{lemma}

\begin{proof}
First we observe that
\begin{align*}
\sum_{\substack{n\leq X\\(n,u)=1}}1
&=\sum_{n\leq X} \sum_{d\mid(n,u)}\mu(d)
=\sum_{\substack{d\mid u}}\mu(d)\sum_{\substack{n\leq X\\ d\mid n}}1
=\sum_{\substack{d\mid u}}\mu(d)\sum_{\substack{n\leq X/d}}1,
\\
\sum_{\substack{n\leq X\\(n,u)=1}}n
&=\sum_{n\leq X} n\sum_{d\mid(n,u)}\mu(d)
=\sum_{\substack{d\mid u}}\mu(d)\sum_{\substack{n\leq X\\ d\mid n}}n
=\sum_{\substack{d\mid u}}\mu(d)\sum_{\substack{n\leq X/d}}nd.
\end{align*}
The expression on the left-hand side of \eqref{lagavulin} is thus equal to
\begin{equation}\label{glenmorangie}
  \sum_{d\mid u}d\mu(d)\left(\frac{X}{d}\sum_{x\leq \frac{X}{d}}1-\sum_{n\leq\frac{X}{d}}n\right).
\end{equation}
Using the identity
\begin{equation*}
    X\sum_{n\leq X}1-\sum_{n\leq X}n=\frac{X^2}{2}-\frac{X}{2}+\frac{\{X\}-\{X\}^2}{2}
=\frac{X^2}{2}-\frac{X}{2}+\frac{\vartheta}{8},
\end{equation*}
we rewrite \eqref{glenmorangie} as follows
\begin{equation*}
    \frac{X^2}{2}\sum_{d\mid u}\frac{\mu(d)}{d}-\frac{X}{2}\sum_{d\mid u}\mu(d)+\frac{\vartheta}{8}\sum_{d\mid u}d\mu^2(d).
\end{equation*}
We conclude the proof by noting that $\sum_{d\mid u}\mu(d) = 0$ since $u>1$.
\end{proof}

\begin{lemma}\label{L:sumstuff}
For a real number $X\geq 1$ and an integer $u>1$, we have
\begin{equation*}
X\sum_{\substack{q\leq X\\u\parallel q}}\frac{\phi(q)}{q}-\sum_{\substack{q\leq X\\u\parallel q}}\phi(q)
=
\frac{3}{\pi^2u}\left(\prod_{p\mid u}\frac{p-1}{p+1}\right)X^2
+ \vartheta \frac{\phi(u) X}{2 u}
\left(
\frac{\sigma(u)}{4} + \frac{\phi(u)}{u} + \frac{\phi(u)}{X}
\right).
\end{equation*}
\end{lemma}

\begin{proof}
We have
\begin{align*}
&\sum_{\substack{q\leq X\\u\parallel q}}\frac{\phi(q)}{q}
=
\sum_{\substack{q\leq \frac{X}{u}\\(q,u)=1}}\frac{\phi(uq)}{uq}
=
\frac{\phi(u)}{u}
\sum_{\substack{q\leq \frac{X}{u}\\(q,u)=1}}\sum_{d\mid q}\frac{\mu(d)}{d}
=
\frac{\phi(u)}{u}
\sum_{\substack{d\leq \frac{X}{u}\\(d,u)=1}}\frac{\mu(d)}{d}\sum_{\substack{q\leq \frac{X}{ud}\\(q,u)=1}} 1,
\\
&\sum_{\substack{q\leq X\\u\parallel q}}\phi(q)
=
\sum_{\substack{q\leq \frac{X}{u}\\(q,u)=1}}\phi(uq)
=
\phi(u)\sum_{\substack{q\leq \frac{X}{u}\\(q,u)=1}}q\sum_{d\mid q}\frac{\mu(d)}{d}
=
\phi(u)
\sum_{\substack{d\leq \frac{X}{u}\\(d,u)=1}}
\mu(d)\sum_{\substack{q\leq \frac{X}{ud}\\(q,u)=1}} q.
\end{align*}
Therefore, using Lemma \ref{bowmore},
\begin{align}\label{glenfiddich}
X\sum_{\substack{q\leq X\\u\parallel q}}\frac{\phi(q)}{q}-\sum_{\substack{q\leq X\\u\parallel q}}\phi(q)
&=
\phi(u)\sum_{\substack{d\leq\frac{X}{u}\\(d,u)=1}}\mu(d)
\left(
\frac{X}{ud}\sum_{\substack{q\leq\frac{X}{ud} \notag \\(q,u)=1}}1
-
\sum_{\substack{q\leq\frac{X}{ud}\\(q,u)=1}}q
\right)
\\
&=
\frac{1}{2}\frac{\phi(u)^2}{u^3}X^2
\sum_{\substack{d\leq \frac{X}{u}\\(d,u)=1}}\frac{\mu(d)}{d^2}
+
\frac{\phi(u)\sigma(u)}{8}\sum_{\substack{d\leq \frac{X}{u}\\(d,u)=1}}\vartheta,
\end{align}
where the absolute value of the second term in \eqref{glenfiddich} is bounded by
\begin{equation*}
    \frac{\phi(u)\sigma(u)}{8u} X.
\end{equation*}
We rewrite the first term in \eqref{glenfiddich} as follows
\begin{align*}
&
\frac{1}{2}\frac{\phi(u)^2}{u^3}X^2
\sum_{\substack{d\leq \frac{X}{u}\\(d,u)=1}}\frac{\mu(d)}{d^2}
\\
\qquad
&=
\frac{1}{2}\frac{\phi(u)^2}{u^3}X^2
\sum_{\substack{(d,u)=1}}\frac{\mu(d)}{d^2}
-
\frac{1}{2}\frac{\phi(u)^2}{u^3}X^2
\sum_{\substack{d>\frac{X}{u}\\(d,u)=1}}\frac{\mu(d)}{d^2}
\\
\qquad
&=
\frac{1}{2}\frac{\phi(u)^2}{u^3}X^2
\frac{6}{\pi^2}\prod_{p\mid u}\left(1-\frac{1}{p^2}\right)^{-1}+\frac{\vartheta}{2}\frac{\phi(u)^2}{u^3}X^2 \left( \frac{u}{X} + \frac{u^2}{X^2} \right)
\\
\qquad
&=
\frac{3}{\pi^2}\frac{1}{u}\left(\prod_{p\mid u}\frac{p-1}{p+1}\right)X^2+\frac{\vartheta}{2}\frac{\phi(u)^2}{u^2}X + \frac{\vartheta}{2} \frac{\varphi(u)^2}{u}.
\end{align*}

\end{proof}
We remark that with more work one could reduce the $\phi(u) \sigma(u)$ term in the conclusion of the lemma by a factor of $6\pi^{-2}$, by using the fact that we are summing over square-free numbers. The bound stated in the lemma is sufficient for our purposes.

We also record the case $u=1$, which is needed in Proposition~\ref{thm:main2} below. The hypothesis $u>1$ in the last two lemmas enters only through the vanishing of $\sum_{d\mid u}\mu(d)$ at the end of the proof of \eqref{lagavulin}; for $u=1$ we have instead
$$X\sum_{n\leq X}1-\sum_{n\leq X}n=\frac{X^2}{2}-\frac{X}{2}+\frac{\vartheta}{8}.$$
Carrying the additional term $-\frac{X}{2d}$ through the proof of Lemma~\ref{L:sumstuff} and using the elementary bound $\bigl|\sum_{d\leq y}\mu(d)/d\bigr|\leq 1$, valid for all $y\geq 1$, shows that the conclusion of Lemma~\ref{L:sumstuff} holds for $u=1$ with the term $\frac{\phi(u)}{u}$ replaced by $2$. This is the source of the convention for $u=1$ in the definition of $E_u(X)$ in Theorem~\ref{T:nonresidues}.
\begin{lemma}\label{prop3}
Let us keep the notations from Theorem \ref{T:nonresidues} and let $N_u(X)$ be as in \eqref{eq: def-N-u}. We have
$$
  N_u(X)\geq E_u(X)\frac{6}{\pi^2}\frac{1}{u}\frac{\phi(u)}{\sigma(u)}X^2h\,.
$$
\end{lemma}
\begin{proof}
Since the number of integer points in $I(q,t) \sqcup J(q,t)$ is at least $2\left(\frac{H}{q} - h\right)$, we have:
\begin{align*}
    N_u(X) &= \sum_{\substack{1 \leq q\leq X\\u\parallel q}} \sum_{\substack{0 \leq t < q\\(t,q)=1}} \sum_{z \in I(q,t) \sqcup J(q,t)} 1 \\
    &\geq 2 \sum_{\substack{1 \leq q\leq X\\u\parallel q}} \left(\frac{H}{q} - h\right) \sum_{\substack{0 \leq t < q\\(t,q)=1}} 1\\
    &= 2h \left( X \sum_{\substack{1 \leq q\leq X\\u\parallel q}} \frac{\varphi(q)}{q} - \sum_{\substack{1 \leq q\leq X\\u\parallel q}} \varphi(q) \right).
\end{align*}
The last lower bound combined with Lemma \ref{L:sumstuff} conclude the proof.
\end{proof}

The following is a modification of \cite[Lemma 5]{Ma}.

\begin{lemma}\label{nicelemma}
Let us keep the notations from Theorem \ref{T:nonresidues} and suppose $3\ell \leq h$. Consider coprime integers $q,t$ satisfying $0 \leq t < q \leq X$. For a fixed choice of cube root of unity $\omega$, assume that $\chi(N)\neq\omega$ for all integers $N\in[1,H]$
with $(N,uv)=1$. If $u\mid q$, then for every integer $z\in\mathcal{I}(q,t)\sqcup\mathcal{J}(q,t)$,
$$
  \left|\sum_{n=0}^{h-1}\chi(z+n)\right|
  \geq
  \frac{h-3\ell}{2}\,.
$$
\end{lemma}

\begin{proof}
We consider only $z\in\mathcal{I}(q,t)$; for $z\in\mathcal{J}(q,t)$ the same argument applies with $q(z+n)-tf$ replaced by its absolute value, which does not change the character value since $\chi(-1)=1$. By Lemma \ref{IJintervals}, if $0\leq n\leq h-1$, then one has $q(z+n)-tf \in (0,H]$. If additionally, $(q(z+n)-tf, uv) = 1$, then $\chi(q(z+n)-tf) \neq \omega$ by assumption. Let us show that $(q(z+n)-tf, uv) \neq 1$ for at most $\ell$ choices of $n$.

If $u_i\mid q(z+n)-tf$ for some $i$, then, since $u_i\mid u\mid q$, we get $u_i\mid tf$; as $(t,q)=1$ we have $u_i\nmid t$, and $u_i\nmid f$ since $f$ is prime and $u_i<h\leq H/2<f$. Hence $(q(z+n)-tf,u)=1$. Thus, it is sufficient to show that for every $1 \leq i \leq \ell$, $v_{i}\mid q(z+n)-tf$ for at most one choice of $n \in [0, h-1]$.

Suppose that $v_{i}\mid q(z+n_1)-tf$ and $v_{i}\mid q(z+n_2)-tf$ for two different values $n_1,n_2 \in [0,h-1]$, then $v_{i}\mid q(n_1 - n_2)$. If $v_i\mid q$, then from $v_i\mid q(z+n_1)-tf$ we get $v_i\mid tf$; since $(t,q)=1$ and $v_i\mid q$, we have $v_i\nmid t$, and hence $v_i\mid f$, contradicting the primality of $f$ (as $v_i<f$). Thus $v_i\nmid q$, and so $v_i$ divides $n_1-n_2$, implying $h\leq v_{i}\leq |n_1-n_2|\leq h-1$, which contradicts the definitions of $n_1$ and $n_2$.

Since $q(z+n)-tf\equiv q(z+n)\pmod{f}$ and $0<q\leq X<f$, we have $\chi(q(z+n)-tf)=\chi(q)\chi(z+n)$, where $\chi(q)$ is a cube root of unity. The argument above thus implies that $\chi(z+n)$, $0 \leq n \leq h-1$, coincides with the cube root of unity $\omega_q:=\chi(q)^{-1}\omega$ for at most $\ell$ values of $n$. Moreover, $\chi(z+n) \neq 0$: writing $q(z+n)=tf+N$ with $0<|N|\leq H<f$, we see that $f\mid z+n$ would force $f\mid N$. Hence we can write
\begin{equation*}
    \sum_{n=0}^{h-1} \chi(z+n) = \omega_q\left(a +b \cdot  \omega +c \cdot  \omega^2\right),
\end{equation*}
where $a$ counts the $n$ with $\chi(z+n)=\omega_q$, so that $a \leq \ell$ and $a+b+c = h$. Since $\omega$ and $\omega^2$ have real part $-\frac12$, we have
\begin{equation*}
\left| \sum_{n=0}^{h-1} \chi(z+n) \right| = \left|a + b\omega + c\omega^2\right| \geq \left|a-\frac{b+c}{2}\right|,
\end{equation*}
where $a - \frac{b+c}{2} = \frac{3a-h}{2} \leq \frac{3\ell-h}{2} \leq 0$ by assumption. This concludes the proof.
\end{proof}

\begin{proof}[Proof of Theorem \ref{T:nonresidues}]
Suppose that $\chi(n)\neq\omega$ for all $n\in[1,H]$ with $(n,uv)=1$.
%
Since $X\geq 2$ and $2HX\leq f$, Lemma~\ref{IJintervals} shows that the intervals $\mathcal{I}(q,t)$, $\mathcal{J}(q,t)$ with $0\leq t<q\leq X$, $(t,q)=1$, are pairwise disjoint subsets of $[-H,f-H)$. The integers $z$ appearing below are therefore pairwise distinct and lie in an interval of length $f$, so by the periodicity of $\chi$, each of them contributes at most one term to $S_\chi(f,h,r)$. Using Lemma \ref{nicelemma} and Lemma \ref{prop3} we find
\begin{align*}
S_\chi(f,h,r)
&=
\sum_{x=0}^{f-1} \left|\sum_{m=0}^{h-1} \chi(x+m) \right|^{2r}
\\
&\geq
\sum_{\substack{0\leq t< q\leq X \\ (q,t) = 1\\u \mid q}}
\sum_{z\in\mathcal{I}(q,t)\sqcup\mathcal{J}(q,t)}
\left|\sum_{m=0}^{h-1}\chi(z+m) \right|^{2r} 
\\
&\geq
\left(\frac{h-3 \ell}{2}\right)^{2r}
N_u(X)
\\
&\geq
E_u(X)
\frac{6}{\pi^2}\frac{1}{u}\frac{\phi(u)}{\sigma(u)}X^2 h
\left(\frac{h-3 \ell}{2}\right)^{2r}
\,.
\end{align*}
Combining this with Proposition \ref{prop2} gives
$$
E_u(X)
\frac{6}{\pi^2}\frac{1}{u}\frac{\phi(u)}{\sigma(u)}X^2 h
\left(\frac{h-3 \ell}{2}\right)^{2r}
\leq
f^{1/2}h^{2r}
W(f,h,r).
$$
We substitute $X=H/h$ and solve for $H^2$ to find
$$
  H^2\leq
  \frac{1}{E_u(X)}
  \frac{\pi^2}{6}u\frac{\sigma(u)}{\phi(u)}
  hf^{1/2}\left(\frac{2h}{h-3 \ell}\right)^{2r}
  W(f,h,r).
$$
\end{proof}

We conclude this section with a complementary bound of P\'olya--Vinogradov strength. It is much weaker than Theorem~\ref{T:nonresidues} when $q_1q_2$ is large, but it applies uniformly when $q_1q_2$ is small, a range in which the hypothesis $2HX\leq f$ prevents us from taking $H$ in Theorem~\ref{T:nonresidues} as large as Criterion~\ref{criterion: expl-non-Eucl} would allow.

\begin{lemma}\label{lem:pv}
Let $\omega$ be a cube root of unity and let $s$ be a squarefree positive integer with $\nu$ prime factors and $(s,f)=1$. If $M$ is a real number satisfying
$$2^{\nu+1}\frac{s}{\phi(s)}\left(f^{1/2}\log f+1\right)\leq M<f,$$
then there exists a positive integer $m\leq M$ with $(m,s)=1$ and $\chi(m)=\omega$.
\end{lemma}
\begin{proof}
Every positive integer $m\leq M$ satisfies $(m,f)=1$, since $M<f$ and $f$ is prime. Hence, writing $N_\omega$ for the number of positive integers $m\leq M$ with $(m,s)=1$ and $\chi(m)=\omega$, orthogonality of the cubic characters and M\"obius inversion give
$$N_\omega=\frac{1}{3}\sum_{j=0}^{2}\bar{\omega}^{j}\sum_{d\mid s}\mu(d)\chi^{j}(d)\sum_{n\leq M/d}\chi^{j}(n).$$
The $j=0$ term equals $\frac{1}{3}\sum_{d\mid s}\mu(d)\lfloor M/d\rfloor\geq\frac{1}{3}\bigl(\frac{\phi(s)}{s}M-2^{\nu}\bigr)$. For $j\in\{1,2\}$, $\chi^{j}$ is a non-principal character modulo $f$, so the P\'olya--Vinogradov inequality (see, e.g., \cite[Ch.~23]{DavMNT}) gives $\bigl|\sum_{n\leq M/d}\chi^{j}(n)\bigr|\leq f^{1/2}\log f$. Therefore
$$N_\omega\geq\frac{1}{3}\left(\frac{\phi(s)}{s}M-2^{\nu}\bigl(2f^{1/2}\log f+1\bigr)\right)>0$$
by the hypothesis on $M$.
\end{proof}
We remark in passing that although numerically sharper versions of the P\'olya--Vinogradov inequality exist (see, e.g., \cite{FS}), what we have above is sufficient for our needs.

\section{The case when \texorpdfstring{$q_1$}{q1} is small}

Here we use the results of the previous section to deal with the case where
the least cubic non-residue $q_1$ is small.  The following result will be 
essential to both our analytic arguments and our computational method.

\begin{theorem} \label{T:1}
For fixed $q_1$, there are no norm-Euclidean cyclic cubic fields for
\mbox{$f\geq f_0(q_1)$}.  For
$q_1\in\{2,3,5,7\}$ one may take $f_0(q_1)=10^{14}$.
See Table~\ref{tab:q1thresh} for admissible values of $f_0(q_1)$ with $q_1\leq 199$.
\end{theorem}
\begin{proof}
Throughout the proof let $Q(f)=1.821 f^{1/4}(\log f)^{3/2}$, so that $q_2\leq Q(f)$ by \cite[Corollary 2]{Ma}, applied with $n = n_0 = 2$ and $p_0 = 10^{10}$.

We first dispose of small $q_2$, uniformly in $q_1\leq 199$. Suppose that
\begin{equation}\label{eq:pvregime}
q_1q_2\leq\frac{f^{1/2}}{73\log f},
\end{equation}
and set $M=24(f^{1/2}\log f+1)$, so that $M<f$ for $f\geq 10^{14}$. Since $s=q_1q_2$ is squarefree with $\nu=2$ prime factors and $\frac{s}{\phi(s)}\leq\frac{2}{1}\cdot\frac{3}{2}=3$, Lemma~\ref{lem:pv} yields a positive integer $m\leq M$ with $(m,q_1q_2)=1$ and $\chi(m)=\chi^{-1}(q_2)$. By \eqref{eq:pvregime},
$$3q_1q_2m\leq\frac{f^{1/2}}{73\log f}\cdot 72\bigl(f^{1/2}\log f+1\bigr)<f
\mand
10q_1^2q_2\leq 10q_1\cdot\frac{f^{1/2}}{73\log f}<f,$$
so Criterion~\ref{criterion: expl-non-Eucl} holds; for $q_1=2$ we use instead $3q_2m\leq 3q_1q_2m<f$.

It remains to treat the range
\begin{equation}\label{eq:bigq2}
\frac{f^{1/2}}{73q_1\log f}<q_2\leq Q(f).
\end{equation}
Here we apply Theorem~\ref{T:nonresidues} with $r=3$, $u=q_1$, $v=q_2$ (thus $\ell=1$), $h=\lceil\lambda f^{1/6}\rceil$ for a suitable $\lambda\in(0.1,2)$ with $h>q_1$, and
\begin{equation}\label{eq:Hmin}
H=\min\left\{\frac{f}{cq_1Q(f)},\,\sqrt{fh/2}\right\},
\end{equation}
where $c=2$ for $q_1>2$, while for $q_1=2$ we take $cq_1=3$ in view of the special shape of Criterion~\ref{criterion: expl-non-Eucl} in that case.
The second term of the minimum guarantees the hypothesis $2HX=2H^2/h\leq f$ of Theorem~\ref{T:nonresidues}. For each $q_1$ and each relevant $\lambda$ we verify that
\begin{equation}\label{eq:seam}
\frac{f^{1/2}}{73q_1\log f}\geq\max\{h,2q_1\},
\end{equation}
so that in the range~\eqref{eq:bigq2} we have $q_2\geq h$, whence $v=q_2$ is admissible in Theorem~\ref{T:nonresidues}, and (for $q_1>2$) $q_2\geq 2q_1$, whence Criterion~\ref{criterion: expl-non-Eucl} requires only $\max\{2q_1q_2m,10q_1^2q_2\}\leq f$ and the first term of the minimum in \eqref{eq:Hmin} gives $H\leq\frac{f}{2q_1q_2}$.

Consider first $q_1=2$, with $\lambda=3/4$ and $f\geq 10^{14}$.
Then $162 \leq h \leq 0.76\,f^{1/6}$, which implies
\begin{equation*}\label{talisker1}
  W(f,h,3)\leq 5+(6+h^{-1})\lambda^{-3}\leq 19.24
\end{equation*}
and
\begin{equation*}\label{talisker2}
\frac{2h}{h-3 \ell} \leq 2.1
\end{equation*}
for $\ell \in \{0, 1, 2\}$. Condition \eqref{eq:seam} reads $\frac{f^{1/2}}{146\log f}\geq h$, which holds for $f\geq 10^{14}$. Both terms of the minimum in \eqref{eq:Hmin} give $X=H/h\geq 1.9\times 10^{5}$: indeed
$$\frac{f/(3Q(f))}{h}\geq\frac{f^{7/12}}{4.16(\log f)^{3/2}}
\mand
\frac{\sqrt{fh/2}}{h}=\sqrt{\frac{f}{2h}}\geq\sqrt{\frac{f^{5/6}}{1.52}},$$
and both right-hand sides exceed $1.9\times 10^{5}$ at $f=10^{14}$ and are increasing in $f$. Hence
$$E_2(X)=1-\frac{\pi^2}{6}\cdot 3\left(\frac{3}{4}+\frac{1}{2}+\frac{1}{X}\right)\frac{1}{X}\geq 0.9995,
\quad\text{so}\quad
\frac{1}{E_2(X)}\frac{\pi^2}{6}\frac{\sigma(2)}{\phi(2)}\cdot 2<10.$$
Therefore the left-hand side of \eqref{eq: THm6-main-condition} is at most
$$10\, h f^{1/2}\times 2.1^6\times 19.24 \cdot H^{-2}<16510\, h f^{1/2} H^{-2}.$$
If $H=f/(3Q(f))$, then $H^2\geq\frac{f^{3/2}}{29.85(\log f)^3}$ and $h\leq 0.76f^{1/6}$, so the right-hand side is less than $1$ provided
$$613\leq\frac{f^{5/12}}{(\log f)^{3/2}},$$
which holds for $f\geq 10^{14}$. If instead $H=\sqrt{fh/2}$, then certainly we have  
$$16510\,hf^{1/2}H^{-2}=33020f^{-1/2}<1.$$ 
Thus Theorem~\ref{T:nonresidues} produces $m\leq H\leq f/(3q_2)$ with $(m,2q_2)=1$ and $\chi(m)=\chi^{-1}(q_2)$, and Criterion~\ref{criterion: expl-non-Eucl} holds. This proves the theorem for $q_1=2$.

The case $q_1=3$ is identical except that $cq_1=6$ in \eqref{eq:Hmin}: condition \eqref{eq:seam} reads $\frac{f^{1/2}}{219\log f}\geq h$, we get $X\geq 9.6\times 10^{4}$ and
$$E_3(X)\geq 1-\frac{11}{X}\geq 0.9995,
\qquad
\frac{1}{E_3(X)}\frac{\pi^2}{6}\frac{\sigma(3)}{\phi(3)}\cdot 3<10,$$
and the first sub-case now requires
$$1225\leq\frac{f^{5/12}}{(\log f)^{3/2}},$$
which again holds for $f\geq 10^{14}$.

Now let $5\leq q_1\leq 199$ and suppose $f\geq f_0\geq 10^{14}$. We proceed computationally. We begin by choosing test points $f_0' \in [10^{14},2\cdot 10^{20}]$. For the purposes of Table~\ref{tab:q1thresh}, these were $f_0' = a \cdot 10^{b}$ where $a \in [1,10)$ in increments of $0.01$. For each test point we varied $\lambda \in (0.1,2)$ and, with $h$ and $H$ as above, checked \eqref{eq:seam} together with the hypotheses and condition \eqref{eq: THm6-main-condition} of Theorem~\ref{T:nonresidues}. In principle, we also need to verify that $10 q_1^2 q_2\leq f$ in the range \eqref{eq:bigq2}. However, using $q_2\leq Q(f)$ allows us to verify this when $f\geq 10^{13}$, for each fixed $3\leq q_1\leq 199$. If these checks pass for a given prime $q_1$, we increment $q_1$ and repeat the checks. We store the largest prime $q_1'$ across all admissible choices for $\lambda \in (0.1,2)$ such that the checks succeed, and define $f_0(q_1)$ to be the smallest test point $f_0'$ such that $q_1\leq q_1'$. A value of $\lambda$ that succeeds at $f_0(q_1)$ is recorded in Table~\ref{tab:q1thresh}. The same verification shows that the checks pass for $q_1\in\{5,7\}$ at every test point down to $10^{14}$. The values of $f_0(q_1)$ could likely be reduced somewhat by a finer search over $\lambda$ and the test points, but this was unnecessary for our application. See \cite{code} for the source code.

It remains to pass from the test points to every $f\geq f_0=f_0(q_1)$, with $q_1$ and $\lambda$ now fixed. Set $\lambda'=\lambda+f_0^{-1/6}$, so that $\lambda f^{1/6}\leq h\leq\lambda'f^{1/6}$ for $f\geq f_0$. Condition~\eqref{eq:seam} follows from $2q_1\leq\lambda'f_0^{1/6}$ together with $f^{1/3}\geq 73q_1\lambda'\log f$, and the latter need only be checked at $f=f_0$ since $f^{1/3}/\log f$ is increasing. For~\eqref{eq: THm6-main-condition}, observe that $h$ is constant on each interval $\bigl(\frac{n-1}{\lambda}\bigr)^{6}<f\leq\bigl(\frac{n}{\lambda}\bigr)^{6}$, and that on such an interval its left-hand side is non-increasing in $f$: the factor $\bigl(\frac{2h}{h-3}\bigr)^{2r}$ is constant, $X=H/h=\min\bigl\{\frac{f}{2q_1Q(f)h},\sqrt{\frac{f}{2h}}\bigr\}$ is increasing so that $E_{q_1}(X)^{-1}$ is decreasing, and
\begin{align*}
\frac{hf^{1/2}W(f,h,3)}{H^{2}}
=\max\Bigl\{\ &4(1.821)^{2}q_1^{2}\Bigl(\frac{5h\log^{3}f}{f}+\frac{(6+h^{-1})\log^{3}f}{h^{2}f^{1/2}}\Bigr),\\
&\frac{10}{f^{1/2}}+\frac{2(6+h^{-1})}{h^{3}}\ \Bigr\},
\end{align*}
each of whose terms is non-increasing for $f\geq e^{6}$. Hence the left-hand side of~\eqref{eq: THm6-main-condition} is bounded on $[f_0,10^{22}]$ by its largest value at $f_0$ and at the left endpoints of the finitely many subsequent intervals; we verify that this largest value is the one at $f_0$. (The left-hand side is not monotone: it can increase by a relative $5\times10^{-5}$ as $h$ increases by $1$.) Finally, replacing $h$ by $\lambda f^{1/6}$ or $\lambda'f^{1/6}$ in the display above according to the sign of its effect yields an explicit decreasing majorant, which does not exceed $6.8\times 10^{-3}$ at $f=10^{22}$; this covers $f\geq 10^{22}$.

\end{proof}


\begin{table}
\begin{tabular}{rcc|rcc|rcc}
$q_1$&$f_0$&$\lambda$&$q_1$&$f_0$&$\lambda$&$q_1$&$f_0$&$\lambda$\\\hline
$2$&$1.00\cdot10^{14}$&$0.75$&$59$&$1.35\cdot10^{17}$&$1.34$&$137$&$3.67\cdot10^{18}$&$1.33$\\
$3$&$1.00\cdot10^{14}$&$0.75$&$61$&$1.54\cdot10^{17}$&$1.34$&$139$&$3.89\cdot10^{18}$&$1.33$\\
$5$&$1.00\cdot10^{14}$&$1.38$&$67$&$2.23\cdot10^{17}$&$1.34$&$149$&$5.11\cdot10^{18}$&$1.33$\\
$7$&$1.00\cdot10^{14}$&$1.38$&$71$&$2.79\cdot10^{17}$&$1.34$&$151$&$5.38\cdot10^{18}$&$1.33$\\
$11$&$2.07\cdot10^{14}$&$1.37$&$73$&$3.11\cdot10^{17}$&$1.34$&$157$&$6.27\cdot10^{18}$&$1.33$\\
$13$&$3.89\cdot10^{14}$&$1.37$&$79$&$4.24\cdot10^{17}$&$1.34$&$163$&$7.27\cdot10^{18}$&$1.33$\\
$17$&$1.08\cdot10^{15}$&$1.36$&$83$&$5.14\cdot10^{17}$&$1.34$&$167$&$7.99\cdot10^{18}$&$1.33$\\
$19$&$1.66\cdot10^{15}$&$1.36$&$89$&$6.76\cdot10^{17}$&$1.34$&$173$&$9.18\cdot10^{18}$&$1.33$\\
$23$&$3.45\cdot10^{15}$&$1.36$&$97$&$9.47\cdot10^{17}$&$1.34$&$179$&$1.05\cdot10^{19}$&$1.33$\\
$29$&$8.47\cdot10^{15}$&$1.35$&$101$&$1.11\cdot10^{18}$&$1.34$&$181$&$1.10\cdot10^{19}$&$1.33$\\
$31$&$1.10\cdot10^{16}$&$1.35$&$103$&$1.20\cdot10^{18}$&$1.33$&$191$&$1.36\cdot10^{19}$&$1.32$\\
$37$&$2.19\cdot10^{16}$&$1.35$&$107$&$1.40\cdot10^{18}$&$1.33$&$193$&$1.42\cdot10^{19}$&$1.32$\\
$41$&$3.26\cdot10^{16}$&$1.35$&$109$&$1.50\cdot10^{18}$&$1.33$&$197$&$1.53\cdot10^{19}$&$1.32$\\
$43$&$3.93\cdot10^{16}$&$1.35$&$113$&$1.73\cdot10^{18}$&$1.33$&$199$&$1.60\cdot10^{19}$&$1.32$\\
$47$&$5.55\cdot10^{16}$&$1.35$&$127$&$2.73\cdot10^{18}$&$1.33$&&&\\
$53$&$8.88\cdot10^{16}$&$1.34$&$131$&$3.08\cdot10^{18}$&$1.33$&&&\\
\end{tabular}
\caption{Values of $f_0(q_1)$ guaranteeing that no norm-Euclidean cyclic cubic field of conductor $f\ge f_0(q_1)$ exists, together with the parameter $\lambda$ used to establish them.}\label{tab:q1thresh}
\end{table}

Using Theorem \ref{T:nonresidues}, we are also able to give a
strong upper bound for $m$ when $q_1$ is large.
In particular, we show the following result, which we will make use of
later.

\begin{prop}
\label{thm:main2}
Let $f\geq 10^{22}$, let $\omega$ be a non-trivial cube root of unity, and suppose $q_1 > 1.82 f^{1/8}$. Then there exists a positive integer $m$ satisfying
$$\chi(m)=\omega, \quad (m,q_1q_2)=1, \quad m\leq 86 f^{5/16}.$$
\end{prop}
\begin{proof}
We apply Theorem \ref{T:nonresidues} with $r=4$, $h = \lceil 1.82 f^{1/8}\rceil$, $H = 86 f^{5/16}$, $u=1$, and $v = q_1q_2$. Since $q_1$ is an integer exceeding $1.82f^{1/8}$, we have $h\leq q_1<q_2$, so $v$ is a product of $\ell=2$ distinct primes in $[h,f)$. Since $f^{1/8}\geq 562$, we have $h\geq 1024$ and $h\leq 1.8218 f^{1/8}$, whence
\begin{equation*}
W(f,h,4) = 7+\frac{f^{1/2}}{h^4}\cdot 24\left(1+\frac{2}{3h}\right)\leq 7+\frac{24}{1.82^4}\left(1+\frac{2}{3072}\right)\leq 9.189
\end{equation*}
and
\begin{equation*}
\left(\frac{2h}{h-6}\right)^{8}\leq\left(\frac{2048}{1018}\right)^{8}\leq 268.4.
\end{equation*}
Moreover, $X=H/h\geq\frac{86}{1.8218}f^{3/16}\geq 6.2\times10^{5}$, so that $E_1(X)\geq 0.999994$, and $2HX=2H^2/h\leq\frac{2\cdot 86^2}{1.82}f^{1/2}<8130 f^{1/2}\leq f$, so that the hypotheses of Theorem~\ref{T:nonresidues} are satisfied.
Hence the left-hand side of \eqref{eq: THm6-main-condition} is at most
\begin{equation*}
\frac{1}{0.999994}\cdot\frac{\pi^2}{6}\cdot\frac{1.8218 f^{5/8}\cdot 268.4\times 9.189}{86^2 f^{5/8}}<0.9994,
\end{equation*}
and the proposition follows.
\end{proof}

\section{Lower bounds for \texorpdfstring{$q_1$}{q1} and \texorpdfstring{$q_2$}{q2}}
\label{sec:mcg}
Having dealt with the cases of $q_1\leq 199$ when $f\geq 10^{22}$
in Theorem~\ref{T:1}, we next focus on proving the following:
\begin{theorem}\label{T:2}
When $q_1>199$, there are no norm-Euclidean cyclic cubic fields with conductor $f \geq 10^{22}$.
\end{theorem}
We assume throughout this section that $q_1\neq 2$ and $f \geq 10^{22}$.
Our proof of Theorem~\ref{T:2} proceeds in a number of stages over the next few sections.
First we give the following lower bounds on $q_1$ and $q_2$.

\begin{prop} \label{fixedprop14}
If the cyclic cubic number field $K$ with conductor $f\geq 10^{22}$ is norm-Euclidean, then $q_1 \ge 379$ and $q_2 \ge 1.3 f^{1/6}$.
\end{prop}

To prove the proposition we will require the following result.
\begin{lemma}
\label{lem:mcgownR}
Let $\omega$ be a primitive cube root of unity. Then for every $k \in \N$, there exists a computable positive constant $D(k)$ such that whenever $f$ satisfies 
\begin{align}
\label{eq:mcgownq}
\frac{4f^{1/4}}{\log^{1/2}f}\geq (2D(k))^k,
\end{align}
there exists a positive integer $m$ satisfying $(m,q_1q_2)=1$, $\chi(m)=\omega$, and 
\begin{equation*}
m<(2D(k))^{k}f^{(k+1)/(4k)}(\log^{1/2} f).
\end{equation*}
If $q_1 \ge 5$ then we may take $D(3)\leq 14.8368$.
\end{lemma}
\begin{proof}
This is a special case of~\cite[Proposition 5.7]{McGown1}.
\end{proof}

\begin{proof}[Proof of Proposition \ref{fixedprop14}]
We proceed via a contrapositive argument.
Assume $q_1 < 379$ or $q_2 < 1.3 f^{1/6}$.
We split the argument into two cases and show in each of them that Criterion \ref{criterion: expl-non-Eucl} holds, and thus $K$ is not norm-Euclidean.
We will assume $q_1\neq 2,3$; indeed the $q_1\in\{2,3\}$ cases of
Theorem~\ref{T:1} permit this whenever $f\geq10^{14}$.
\vspace{1ex}



\noindent
{\bf Case I:}  $3 < q_1 < q_2 < 1.3 f^{1/6}$.
The hypothesis \eqref{eq:mcgownq} holds for $k=3$ since $4f^{1/4}/\log^{1/2}f\geq 1.7\times10^{5}>(2\cdot14.8368)^3$ for $f\geq10^{22}$. Thus by Lemma \ref{lem:mcgownR} applied with $\omega=\chi^{-1}(q_2)$, for $5 \le q_1 < q_2$, there exists $m$ with $(m,q_1q_2)=1$, $\chi(m)=\omega$, and
\begin{equation*}
m < (2 \cdot 14.8368)^3 f^{1/3} (\log^{1/2} f) < 26129 f^{1/3} (\log^{1/2} f) \leq \frac{f}{3 q_1 q_2},
\end{equation*}
where the last bound is implied by $q_1 < q_2 < 1.3 f^{1/6}$ and $f \geq 10^{22}$. Moreover, $10q_1^2q_2 < 10\cdot(1.3)^3f^{1/2}\leq f$. We conclude that $K$ is not norm-Euclidean by Criterion \ref{criterion: expl-non-Eucl}.

\vspace{1ex}
\noindent
{\bf Case II:}  $q_1 < 379 < 1.3 f^{1/6} \le q_2$.
Then $q_1 \leq 373$ since it is prime and $q_2 \ge 1.3 f^{1/6} \ge 6034 > 2 q_1$. Recalling from the proof of Theorem~\ref{T:1} that $q_2\leq Q(f):=1.821f^{1/4}\log^{3/2}f$, we have $10q_1^2q_2\leq 10\cdot 373^2\,Q(f)\leq f$ for $f\geq 10^{22}$. For each prime $q_1 \le 373$, we invoke Theorem \ref{T:nonresidues} with $r = 3$,
$$H = \min\left\{\frac{f}{2 \times  373 \times Q(f)},\,\sqrt{fh/2}\right\} \leq \frac{f}{2 q_1 q_2}\,,$$
$u=q_1$, $v=q_2$ (so that $\ell=1$), $h = \lceil 1.3 f^{1/6} \rceil$ and $\omega=\chi^{-1}(q_2)$; note that $v=q_2$ is admissible since $q_2$ is an integer with $q_2\geq 1.3f^{1/6}$, whence $q_2\geq\lceil 1.3f^{1/6}\rceil=h$. The second term of the minimum ensures the hypothesis $2HX\leq f$ of Theorem \ref{T:nonresidues}; it is the smaller of the two only for $f>10^{35}$, where the left-hand side of \eqref{eq: THm6-main-condition} is below $10^{-11}$.
We verify all of the hypotheses of Theorem \ref{T:nonresidues} computationally \cite{code} throughout the range: one has $X\geq 1.07\times10^{7}$, $X/q_1\geq 2.8\times10^{4}$ and $E_{q_1}(X)\geq 0.994$, while $h\geq 6035$ makes $r\leq 9h$ and $3\ell<h$ immediate. The left-hand side of \eqref{eq: THm6-main-condition} never exceeds $0.045$, thus implying Criterion \ref{criterion: expl-non-Eucl}.
\end{proof}
\section{Character sum estimates}
\label{charsum}
In this section we establish a new character sum estimate (see Theorem~\ref{thm:1}), which will be combined with Vinogradov's sieve in Section~\ref{sec:vino} to obtain improved bounds for $q_1,q_2$ in the cases not covered by Theorem~\ref{T:1} and Proposition~\ref{fixedprop14}. Finally, in Section~\ref{sec:II} we will complete the proof of Theorem~\ref{T:2}. The main result of this section is as follows.
\begin{theorem}
\label{thm:1}
Let $f$ be prime and $U,V,T,r \in \N$. Consider the character $\chi$ modulo $f$ of order $3$ and some real number $0<\theta<1$. Define the quantities
\begin{align}
&S=\sum_{1\leq t \leq T}\sum_{1\leq u \leq U}\sum_{1\leq v \leq V}\chi(v-ut),\nonumber \\
&\Delta=\left(\frac{f^{1/2}W(f,T,r)}{UV}\right)^{1/2r} \label{eq:Deltadef}, \\
&
E(U, V)= \frac{6}{\pi^2} \left( U^{-1} + V^{-1} \right) \log^2 U + \left( 2.66 U^{-1} + 2.26 V^{-1} \right) \log U \label{eq:EUt}\\
&\quad \quad \quad \,\, \,+ 7.57 U^{-1} + 9.82 V^{-1},\nonumber\\
& E_1 = \Delta \times  \left(\frac{12}{\pi^2} \theta \log U + 1 + U^{\theta} \times  E(U,V)\right)^{1/2r},\nonumber\\
& E_2 = \frac{12}{\pi^2} \, U^{-\theta} + E(U,V).\nonumber
\end{align}
Suppose $UV < f$, then $|S|\leq TUV (E_1 + E_2)$.
\end{theorem}
We provide some preliminary estimates required for the proof of Theorem \ref{thm:1}. What follows is essentially due to Bourgain, Garaev, Konyagin and Shparlinski~\cite[Lemma 5]{BGKS}. Specialising to intervals starting at the origin allows for sharper numerical results. 
\begin{lemma}
\label{lem:1}
Let $f$ be prime and $U,V \in \N$ satisfying $UV < f$. Let $\lambda \not \equiv 0 \dmod{f}$ and for $1\leq \Delta$, suppose the congruence relation
\begin{equation} \label{eq:lem1cong}
v\equiv \lambda u \dmod{f}, \quad 1 \leq u \leq U, \quad 1 \leq v \leq V,
\end{equation}
has at least $\Delta$ solutions. Then there exist $u_0, v_0 \in \N$ satisfying 
\begin{align}
\label{eq:n0u0bound}
1 \leq u_0\leq \frac{U}{\Delta}, \quad 1 \leq v_0\leq \frac{V}{\Delta},  \quad (u_0, v_0)=1,
\end{align}
such that 
\begin{align*}
v_0 \equiv \lambda u_0 \dmod{f}.
\end{align*}
\end{lemma}
\begin{proof}
The congruence \eqref{eq:lem1cong} has at least one solution $(u, v)$ since $1 \leq \Delta$. Let
\begin{equation*}
    u_0 = \frac{u}{\gcd(u,v)}, \quad v_0 = \frac{v}{\gcd(u,v)}.
\end{equation*}
Since $1 \leq uv < f$, the numbers $f$ and $\gcd(u, v)$ are coprime, whence $(u_0, v_0)$ is also a solution to \eqref{eq:lem1cong}. We will show $u_0, v_0$ satisfy the inequalities in \eqref{eq:n0u0bound}. By the assumptions in the lemma,
\begin{align*}
|\{ (u, v)\in [1,U]\times [1,V] \ : \  u_0 v\equiv u v_0 \dmod{f}\}|\geq \Delta,
\end{align*}
which is equivalent to
\begin{equation} \label{eq1: lem BGKS}
|\{ (u, v)\in [1,U]\times [1,V] \ : \  u_0 v=u v_0\}|\geq \Delta,
\end{equation}
since $1 \leq UV < f$. Let $(u, v)$ be from the set in \eqref{eq1: lem BGKS}. Since $(u_0, v_0) = 1$, there exists $\ell \in \N$ such that $u = \ell u_0$ and $v = \ell v_0$. Hence, \eqref{eq1: lem BGKS} can be rewritten as follows
\begin{align*}
|\{ \ell \in \N \ : \  1\leq \ell u_0\leq U \quad \text{and} \quad 1\leq \ell v_0\leq V \}|\geq \Delta.
\end{align*}
In particular, this implies two inequalities
\begin{equation*}
|\{ \ell \in \N \ : \  1\leq \ell u_0\leq U  \}|\geq \Delta, \quad |\{ \ell \in \N \ : \  1\leq \ell v_0\leq V  \}|\geq \Delta,
\end{equation*}
and thus
\begin{equation*}
\Delta \leq \frac{U}{u_0}, \quad \Delta \leq \frac{V}{v_0},
\end{equation*}
which completes the proof of the lemma.
\end{proof}
The next lemma is elementary, though we are not aware of it appearing in the literature before now.
\begin{lemma}
\label{lem:coprimepair}
For any positive integers $U$ and $V$ we have 
\begin{equation} \label{eq1: lem coprime pairs}
|\{ (u, v)\in [1,U]\times [1,V] \ : \ (u,v)=1\}| \leq B_1(U, V),
\end{equation}
where
\begin{equation} \label{eq2: lem coprime pairs def B1(U,V)}
B_1(U, V) := \frac{6}{\pi^2} UV + \frac{6}{\pi^2} (U+V) \log U + 2.044 V + 1.652 U + 2.04 \frac{V}{\sqrt{U}} + 2.72 \sqrt{U}.
\end{equation}
\end{lemma}
\begin{proof}
Let $S$ denote the set on the left-hand side of \eqref{eq1: lem coprime pairs}. Let $\gamma = 0.577\ldots$ denote Euler's constant. Then
\begin{align*}
S&=\sum_{1\leq u \leq U}\sum_{\substack{1 \leq v\leq V \\ (u, v)=1}}1
= \sum_{1\leq u \leq U} \sum_{1 \leq v\leq V} \sum_{d \mid (u, v)} \mu(d) = \sum_{1 \leq d \leq U} \mu(d) \left\lfloor \frac{U}{d} \right\rfloor \left\lfloor \frac{V}{d} \right\rfloor \\
&\leq UV \sum_{1 \leq d \leq U} \frac{\mu(d)}{d^2} + (U+V) \sum_{1 \leq d \leq U} \frac{\mu^2(d)}{d} + \sum_{1 \leq d \leq U} \mu^2(d) \\
&\leq UV \left( \frac{1}{\zeta(2)} + \frac{1}{U} \right) + (U + V) \left( \frac{\log U}{\zeta(2)} + \frac{\gamma}{\zeta(2)} - 2 \frac{\zeta'(2)}{\zeta^2(2)} + \frac{2.04}{\sqrt{U}} \right) + \frac{U}{\zeta(2)} + 0.68 \sqrt{U},
\end{align*}
where we used $\sum_{d \geq 1} \frac{\mu(d)}{d^2} = \frac{1}{\zeta(2)}$, \cite[(4.11)]{Buthe1}, and \cite[(4.6)]{Buthe1} respectively.
We conclude by noting that $\zeta(2) = \frac{\pi^2}{6}$ and $\zeta'(2) \geq -0.9376$.
\end{proof}
\begin{rem}\label{rem:B1real}
The bound \eqref{eq1: lem coprime pairs} remains valid for real $U,V\geq1$: its left-hand side is unchanged upon replacing $U,V$ by $\lfloor U\rfloor,\lfloor V\rfloor$, and $B_1$ is non-decreasing in each argument, since every term is non-decreasing in $V$, while the only term decreasing in $U$ satisfies
$$\left|\frac{\partial}{\partial U}\frac{2.04V}{\sqrt{U}}\right|=\frac{1.02V}{U^{3/2}}\leq\frac{6}{\pi^2}V\left(1+\frac1U\right)\leq\frac{\partial}{\partial U}\left(\frac{6}{\pi^2}UV+\frac{6}{\pi^2}(U+V)\log{U}\right)$$
for $U\geq1$. For $0<V<1$ the left-hand side of \eqref{eq1: lem coprime pairs} vanishes, so the bound holds trivially. We use Lemma \ref{lem:coprimepair} below with real arguments in this sense.
\end{rem}

\begin{proof}[Proof of Theorem \ref{thm:1}]
Since $1 \leq u \leq U < f$, $u$ is invertible modulo $f$. Let $\overline{u}$ denote an inverse to $u$ modulo $f$, then
\begin{align*}
|S|&\leq \sum_{\substack{1\leq u \leq U \\ 1\leq v \leq V}}\left|\sum_{1\leq t \leq T}\chi(\overline{u}v-t)\right| \\
&\leq \sum_{\lambda=1}^{f-1}I(\lambda)\left|\sum_{1\leq t \leq T}\chi(\lambda - t)\right|,
\end{align*}
where 
\begin{align*}
I(\lambda)=\#\{ 1\leq u \leq U, \ \ 1\leq v \leq V \ : \ v\equiv \lambda u \dmod{f}\}.
\end{align*} 
Trivially, we have $I(\lambda)\leq U$. For each $j$  define the set 
\begin{equation*}
D_j=\{ 1\leq \lambda \leq f-1 \ : I(\lambda)=j\}.
\end{equation*}
We partition $S$ with respect to the parameter $\theta$ as follows
\begin{align}
\label{eq:S12}
|S|\leq \sum_{\substack{1\leq j \leq U^{\theta} }}\sum_{\lambda \in D_j}I(\lambda)\left|\sum_{1\leq t \leq T}\chi(\lambda - t)\right|+\sum_{\substack{ j> U^{\theta} }}\sum_{\lambda \in D_j}I(\lambda)\left|\sum_{1\leq t \leq T}\chi(\lambda - t)\right| =:  S_1+S_2.
\end{align}

\noindent
{\bf Estimation of $S_2$.}
\begin{align*}
S_2&\leq T\sum_{\substack{U^{\theta}\leq j \leq U}}j|D_j|=T\sum_{0\leq j \leq U-U^{\theta}}(U-j)|D_{U-j}| \\
&=TU\sum_{0\leq j \leq U-U^{\theta}}|D_{U-j}|-T\sum_{0\leq j \leq U-U^{\theta}}j|D_{U-j}|.
\end{align*}
By partial summation 
\begin{align*}
\sum_{0\leq j \leq U-U^{\theta}}j|D_{U-j}|=(U-U^{\theta})\sum_{0\leq j \leq U-U^{\theta}}|D_{U-j}|-\int_{0}^{U-U^{\theta}}\sum_{0\leq j \leq w}|D_{U-j}|dw,
\end{align*}
and hence 
\begin{align}
\label{eq:S2b111}
S_2\leq T \times  g(U, U^{\theta}),
\end{align}
where 
\begin{equation*} 
g(U, Y) = Y\sum_{0\leq j \leq U-Y}|D_{U-j}|+\int_{0}^{U-Y}\sum_{0\leq j \leq w}|D_{U-j}|dw.
\end{equation*}
Note that $I(\lambda)\geq (U-w)$ is equivalent to $\lambda \in D_{U-j}$ for some $0\leq j \leq w$, hence
\begin{align*}
\sum_{0\leq j \leq w}|D_{U-j}| = |\{ 1 \leq \lambda \leq f-1 \ :  \ I(\lambda)\geq U-w\}|.
\end{align*}
Let $D_{\geq U - w}$ denote the set on the right-hand side of the equation above. For each $\lambda \in D_{\geq U - w}$, we apply Lemma \ref{lem:1} with $\Delta=U-w$ to obtain a pair $(u_\lambda,v_\lambda)$ of positive integers satisfying 
\begin{equation*}
v_\lambda \equiv \lambda u_\lambda \dmod{f}, \quad 1 \leq u_\lambda\leq \frac{U}{(U-w)}, \quad 1\leq v_\lambda\leq \frac{V}{(U-w)}, \quad (u_\lambda,v_\lambda)=1.
\end{equation*}
Since for $\lambda_1 \neq \lambda_2$, the pairs $(u_{\lambda_1},v_{\lambda_1}) \neq (u_{\lambda_2},v_{\lambda_2})$, we have:
\begin{equation*}
    |D_{\geq U - w}| \leq \left|\left\{ (u, v)\in \left[1,\frac{U}{U-w}\right]\times \left[1,\frac{V}{U-w}\right] \ : \ (u,v)=1\right\}\right| \leq B_1\left( \frac{U}{U-w}, \frac{V}{U-w} \right),
\end{equation*}
with $B_1$ defined in \eqref{eq2: lem coprime pairs def B1(U,V)}. Hence by \eqref{eq:S2b111},
\begin{align*}
S_2\leq T U^{\theta} \times  B_1(U^{1-\theta}, U^{-\theta}V) + T \int_{0}^{U-U^{\theta}} B_1\left( \frac{U}{U-w}, \frac{V}{U-w} \right) dw.
\end{align*}
Let us bound the last integral above: for $1 \leq Y \leq U$, we have
\begin{align*}
    \int_{0}^{U-Y}& B_1\left( \frac{U}{U-w}, \frac{V}{U-w} \right) dw \\
    &= \frac{6}{\pi^2} U V Y^{-1} - \log Y \left( \frac{6}{\pi^2} (U+V)\log U + 1.652 U + 2.044 V\right) \\
    &+ \frac{3}{\pi^2} (U+V) \log^2 Y - 2 Y^{1/2} \left( 2.04 \frac{V}{\sqrt{U}} + 2.72 \sqrt{U} \right) \\
    &+ \log U (1.652 U + 2.044 V) + \frac{3}{\pi^2} (U+V) \log^2 U + 5.44 U + V \left( 4.08 - \frac{6}{\pi^2} \right) \\
    &\leq \frac{6}{\pi^2} UVY^{-1} + \frac{6}{\pi^2} (U+V) \log^2 U + (1.652 U + 2.044 V) \log U + 5.44 U + 3.48 V.
\end{align*}
In addition, since $1\leq Y\leq U$,
\begin{align*}
    Y\, B_1\left( \frac{U}{Y}, \frac{V}{Y}\right) &= \frac{6}{\pi^2} UVY^{-1} + \frac{6}{\pi^2} (U+V) \log\frac{U}{Y} + 1.652U + 2.044V \\
    &\quad + 2.04 V\sqrt{\frac{Y}{U}} + 2.72\sqrt{UY}\\
    &\leq \frac{6}{\pi^2} UVY^{-1} + \frac{6}{\pi^2} (U+V) \log U + 4.372 U + 4.084 V,
\end{align*}
whence, taking $Y=U^{\theta}$ and absorbing $\frac{6}{\pi^2}(U+V)\log U$ into the coefficients of $U\log U$ and $V\log U$ (note $1.652+\frac{6}{\pi^2}<2.26$ and $2.044+\frac{6}{\pi^2}<2.66$), we get
\begin{align} \label{eq:S2b11}
    S_2 \leq TUV &\bigg( \frac{12}{\pi^2} U^{-\theta} + \frac{6}{\pi^2} \left( U^{-1} + V^{-1} \right) \log^2 U \\
    &\,\, + \left( 2.66 U^{-1} + 2.26 V^{-1} \right) \log U + 7.57 U^{-1} + 9.82 V^{-1} \bigg). \nonumber
\end{align}

\noindent
{\bf Estimation of $S_1$.} Let $p = \frac{2r}{2r-1}$, $q = 2r$, and
\begin{equation*}
    a_\lambda = \left| \sum_{1 \leq t \leq T} \chi(\lambda - t) \right|,
\end{equation*}
then
\begin{equation*}
    S_1 = \sum_{1 \leq j \leq U^{\theta}} j \sum_{\lambda \in D_j} a_\lambda.
\end{equation*}
By H\"{o}lder's inequality,
\begin{equation*}
    \sum_{\lambda \in D_j} a_\lambda \leq |D_j|^{1/p} \left(\sum_{\lambda \in D_j} a_\lambda^{q}\right)^{1/q},
\end{equation*}
whence, by applying H\"{o}lder's inequality two more times, we get
\begin{align*} 
S_1 &\leq \sum_{1 \leq j \leq U^{\theta}} j |D_j|^{1/p} \left(\sum_{\lambda \in D_j} a_\lambda^{q}\right)^{1/q} \\
&\leq \left( \sum_{1 \leq j \leq U^{\theta}} j^p |D_j| \right)^{1/p} \times  \left( \sum_{1 \leq j \leq U^{\theta}} \left( \sum_{\lambda \in D_j} a_\lambda^{q} \right) \right)^{1/q} \\
&= \left( \sum_{1 \leq j \leq U^{\theta}} (j |D_j|)^{2-p} \times  (j^2 |D_j|)^{p-1} \right)^{1/p} \times  \left( \sum_{1 \leq j \leq U^{\theta}} \left( \sum_{\lambda \in D_j} a_\lambda^{q} \right) \right)^{1/q}\\
&\leq \left( \sum_{1 \leq j \leq U^{\theta}} j |D_j| \right)^{(2-p)/p} \times  \left(\sum_{1 \leq j \leq U^{\theta}} j^2 |D_j| \right)^{(p-1)/p} \times  \left( \sum_{1 \leq j \leq U^{\theta}} \left( \sum_{\lambda \in D_j} a_\lambda^{q} \right) \right)^{1/q}.
\end{align*}
The last chain of inequalities implies
\begin{equation} \label{eq: bound-S1-Holder}
    (S_1)^{2r} \leq \left(\sum_{1 \leq j \leq U^{\theta}}j |D_j|\right)^{2r-2}\left(\sum_{1 \leq j \leq U^{\theta}}j^2|D_j|\right) \left(\sum_{\lambda=1}^{f}\left|\sum_{1\leq t \leq T}\chi(\lambda-t)\right|^{2r}\right).
\end{equation}
We can bound the last factor on the right-hand side of \eqref{eq: bound-S1-Holder} by Proposition \ref{prop2},
\begin{equation} \label{eq: bound-S1-term3}
\sum_{\lambda=1}^{f}\left|\sum_{1\leq t \leq T}\chi(\lambda - t)\right|^{2r} = \sum_{\lambda=1}^{f}\left|\sum_{0\leq t \leq T-1}\chi(\lambda + t)\right|^{2r}
\leq f^{1/2}T^{2r}W(f,T,r).
\end{equation}
We bound the first factor in \eqref{eq: bound-S1-Holder} using
\begin{equation} \label{eq: bound-S1-term1}
\sum_{1 \leq j \leq U^{\theta}}j |D_{j}|\leq \sum_{1\leq \lambda \leq f}I(\lambda) \leq UV,
\end{equation}
since the solution sets counted by $I(\lambda)$, $1 \leq \lambda \leq f-1$, are pairwise disjoint subsets of $[1,U] \times [1,V]$. Finally, let us bound the second factor on the right-hand side of \eqref{eq: bound-S1-Holder}. By partial summation,
\begin{align*}
\sum_{1 \leq j \leq U^{\theta}}j^2|D_{j}|=U^{\theta}\sum_{1 \leq j \leq U^{\theta}}j|D_{j}|-\int_{1}^{U^{\theta}}\sum_{1\leq j\leq w}j|D_{j}|dw,
\end{align*} 
and
\begin{align*}
\sum_{1 \leq j \leq w}j |D_{j}|=UV-\sum_{w < j \leq U}j |D_{j}|,
\end{align*}
so that 
\begin{align} \label{eq: bound-j^2-Dj}
\sum_{1 \leq j \leq U^{\theta}}j^2|D_{j}|&=UV + \int_{1}^{U^{\theta}}\sum_{w < j \leq U}j |D_{j}|dw - U^{\theta}\sum_{U^{\theta} < j \leq U}j |D_{j}| \leq UV + \int_{1}^{U^{\theta}}\sum_{w \leq j \leq U}j |D_{j}|dw.
\end{align}

Arguing similarly to the proof of \eqref{eq:S2b11}, for $w\geq 1$ we have
\begin{align*}
\sum_{w \leq j \leq U}j|D_{j}|\leq UV &\bigg( \frac{12}{\pi^2} w^{-1} + \frac{6}{\pi^2} \left( U^{-1} + V^{-1} \right) \log^2 U \\
    &\,\, + \left( 2.66 U^{-1} + 2.26 V^{-1} \right) \log U + 7.57 U^{-1} + 9.82 V^{-1} \bigg),
\end{align*}
and therefore 
\begin{align*}
\int_{1}^{U^{\theta}}\sum_{w \leq j \leq U}j |D_{j}|dw \leq 
UV &\bigg( \frac{12}{\pi^2} \theta \log U + \frac{6}{\pi^2} \left( U^{-1} + V^{-1} \right) U^{\theta}\log^2 U \\
    &\,\, + \left( 2.66 U^{-1} + 2.26 V^{-1} \right) U^{\theta} \log U + 7.57 U^{\theta -1} + 9.82 U^{\theta}V^{-1} \bigg).
\end{align*}
Combining with \eqref{eq: bound-j^2-Dj}, we obtain
\begin{equation} \label{eq: bound-S1-term2}
\sum_{1 \leq j \leq U^{\theta}}j^2|D_{j}| \leq 
UV \left( \frac{12}{\pi^2} \theta \log U + 1 + U^{\theta} \times  E(U,V) \right),
\end{equation}
where $E(U,V)$ is defined in \eqref{eq:EUt}.

The bounds \eqref{eq: bound-S1-Holder}, \eqref{eq: bound-S1-term3}, \eqref{eq: bound-S1-term1}, and \eqref{eq: bound-S1-term2} imply
\begin{equation*}
S_1 \leq TUV \times  \Delta \times \left( \frac{12}{\pi^2} \theta \log U + 1 + U^{\theta} \times  E(U,V) \right)^{1/2r},
\end{equation*}
with $\Delta$ defined in \eqref{eq:Deltadef}. Combining this with \eqref{eq:S12} and \eqref{eq:S2b11} completes the proof.
\end{proof}

\section{Cubic non-residues II}
\label{sec:vino}
The aim of this section is to improve on the upper bounds for $q_1,q_2$ whenever $q_1$ is bounded uniformly from below. Our approach is to combine the character sum estimates from Section~\ref{charsum} with Vinogradov's sieve. We recall that by 
\cite{Lezowski},
all cyclic cubic norm-Euclidean number fields have conductor
$f \leq 10^{50}$,
and therefore we may assume this upper bound on $f$.
\begin{theorem}
\label{thm:main1}
Let $10^{22} \leq f \leq 10^{50}$ be prime.
If $q_1\geq 379$, then
$$q_{1} \leq q_{2} \leq 36.88  f^{0.21769}.$$
\end{theorem}
The proof of Theorem~\ref{thm:main1} requires some preliminary results.
\begin{lemma}
\label{lem:dus1}
Let 
\begin{align*}
B=\gamma+\sum_{p}\left(\log\left(1-\frac{1}{p} \right)+\frac{1}{p} \right)= 0.26149\ldots.
\end{align*}
Then
\begin{align*}
\left|\sum_{p\leq x}\frac{1}{p} -\log\log{x}-B\right| \leq \frac{1}{2\log^{2}{x}}\quad \text{for } x \geq 286,
\end{align*}
and moreover
\begin{align*}
\sum_{p\leq x}\frac{1}{p} > \log\log{x}+B \quad \text{for } 2 < x \leq 10^{12}.
\end{align*}
\end{lemma}
\begin{proof}
See \cite[Theorem~5]{R-S} and \cite[Lemma~2.13]{BJS22}.
\end{proof}
\begin{lemma}
\label{lem:firsttwo11}
Let $T,U,V$ be positive integers, and let $V_0 = \max\{V, UT\}$; suppose $V_0<f$. For $x>1$ let $\delta_0(x) = 0$ if $x\leq 10^{12}$ and $\delta_0(x)=\frac{1}{2\log^2 x}$ otherwise. Let $\varepsilon>0$ and suppose
\begin{align}
\label{eq:firsttwoassumption}
\left|\sum_{\substack{1\leq t \leq T \\ 1\leq u \leq U \\ 1\leq v \leq V}}\chi(v-ut)\right|\leq \varepsilon TUV.
\end{align}
If $q_2 \geq 286$ then we have 
\begin{align}
\label{eq:fft1}
\frac{1}{q_1}+\frac{1}{q_2}+\log\left(\frac{\log{V_0}}{\log{q_2}}\right)+\delta_0(q_2) + \frac{1}{2\log^2{V_0}} + \frac{\pi(V_0)}{V} \geq \frac{2(1-\varepsilon)}{3},
\end{align}
and if moreover $q_1\geq 286$ then we have 
\begin{align}
\label{eq:fft2}
\frac{1}{q_1} + \log\left(\frac{\log{V_0}}{\log{q_1}}\right) + \delta_0(q_1) + \frac{1}{2\log^2{V_0}} + \frac{\pi(V_0)}{V}\geq \frac{2(1-\varepsilon)}{3}.
\end{align}
\end{lemma}
\begin{proof}
The proof is based on Vinogradov's trick~\cite{Vin}. Let us define
\begin{align*}
&A_0=\{ (t,u,v)\in [1,T]\times [1,U]\times [1,V] \ : \ \chi(v-ut) = 1 \},\\
&A_0^{c}= [1,T]\times [1,U]\times [1,V] \setminus A_0.
\end{align*}
Let $\chi_0$ denote the principal character modulo $f$, and $\overline{\chi}$ the complex-conjugate character of $\chi$. Then we can express $A_0$ as follows
\begin{align*}
|A_0|=\frac{1}{3}\sum_{\psi \in \{\chi, \overline{\chi}, \chi_0\}}\sum_{\substack{1\leq t \leq T \\ 1\leq u \leq U \\ 1\leq v \leq V}}\psi(v-ut).
\end{align*}
By isolating the contribution from $\chi_0$ and using the assumption \eqref{eq:firsttwoassumption} for $\chi, \overline{\chi}$, we get 
\begin{align*}
|A_0| \leq \frac{1+2\varepsilon}{3}TUV,
\end{align*}
and hence 
\begin{align}
\label{eq:T3UB}
|A_0^{c}|\geq \frac{2(1-\varepsilon)}{3}TUV.
\end{align}
We note that $|v-ut| \leq V_0$ for all $(u,v,t) \in A_0^{c}$. If $\chi(v-ut)=0$ then $f\mid v-ut$, which forces $v=ut$ since $|v-ut|\leq V_0<f$; such triples are counted by the first sum below. Otherwise $v-ut$ is a cubic non-residue, and by the definitions of $q_1, q_2$ every cubic non-residue not exceeding $V_0$ in absolute value must be divisible by a prime $p$ satisfying $p=q_1$ or $q_2\leq p \leq V_0$. This implies
\begin{align}
\label{eq:T3lb1}
|A_0^{c}|\leq \sum_{\substack{1\leq t \leq T \\ 1\leq u \leq U \\ 1\leq v \leq V \\ v-ut\equiv 0 \dmod{q_1}}}1+ \sum_{q_2\leq p \leq V_0}\sum_{\substack{1\leq t \leq T \\ 1\leq u \leq U \\ 1\leq v \leq V \\ v-ut\equiv 0 \dmod{p}}}1,
\end{align}
and a slightly weaker bound
\begin{align}
\label{eq:T3lb2}
|A_0^{c}|\leq \sum_{q_1\leq p \leq V_0}\sum_{\substack{1\leq t \leq T \\ 1\leq u \leq U \\ 1\leq v \leq V \\ v-ut\equiv 0 \dmod{p}}}1.
\end{align}
For any prime $p$ and any pair $(t,u)$, the number of $v\in[1,V]$ with $v\equiv ut\dmod{p}$ is at most $V/p+1$. Hence the bound \eqref{eq:T3lb1} implies
\begin{align*}
|A_0^{c}|&\leq TUV\left(\frac{1}{q_1}+\sum_{q_2\leq p \leq V_0}\frac{1}{p}\right) + TU\bigl(\pi(V_0)-\pi(q_2)+2\bigr)\\
&\leq TUV\left(\frac{1}{q_1}+\sum_{q_2\leq p \leq V_0}\frac{1}{p} + \frac{\pi(V_0)}{V} \right),
\end{align*}
using $\pi(q_2)\geq 2$ in the last step. Writing
\begin{align*}
\sum_{q_2\leq p \leq V_0}\frac{1}{p} = \sum_{p \leq V_0}\frac{1}{p}-\sum_{p \leq q_2}\frac{1}{p}+\frac{1}{q_2}
\end{align*}
and applying both parts of Lemma \ref{lem:dus1} (according to whether $q_2\leq 10^{12}$), we get for $q_2 \geq 286$
\begin{align*}
|A_0^{c}|\leq TUV\left(\frac{1}{q_1}+\frac{1}{q_2}+\log\left(\frac{\log{V_0}}{\log{q_2}}\right)+\delta_0(q_2) + \frac{1}{2\log^2{V_0}} + \frac{\pi(V_0)}{V}\right),
\end{align*}
which combined with \eqref{eq:T3UB} implies \eqref{eq:fft1}. Using \eqref{eq:T3lb2} and a similar argument with $q_1 \geq 286$ we get \eqref{eq:fft2}.
\end{proof}
\begin{cor}
\label{lem:firsttwo1}
Keep the notations and assumptions from Lemma \ref{lem:firsttwo11}. Let us define $\rho_1$ and $\rho_2$ implicitly by
\begin{align*}
q_1=V_0^{\rho_1}, \quad q_2=V_0^{\rho_2}.
\end{align*}
Then if $q_2 \geq 286$,
\begin{align*}
\log\left(\frac{1}{\rho_2}\right)\geq \frac{2(1-\varepsilon)}{3}-\frac{1}{q_1}-\frac{1}{q_2}-\delta_0(q_2)-\frac{1}{2\log^{2} V_0} - \frac{\pi(V_0)}{V},
\end{align*}
and if $q_1 \geq 286$,
\begin{align*}
\log\left(\frac{1}{\rho_1}\right)\geq \frac{2(1-\varepsilon)}{3} - \frac{1}{q_1}-\delta_0(q_1)-\frac{1}{2\log^{2} V_0} - \frac{\pi(V_0)}{V}.
\end{align*}
\end{cor}

\begin{proof}[Proof of Theorem \ref{thm:main1}]
By Theorem \ref{thm:1}, we have 
\begin{equation*}
\left|\sum_{\substack{1\leq t \leq T \\ 1\leq u \leq U \\ 1\leq v \leq V}}\chi(v-ut)\right|\leq TUV(E_1+E_2),
\end{equation*}
where we choose parameters
\begin{align*}
    r = 2, \quad \theta = 0.44, \quad \alpha=\sqrt{3/2}, \quad K_1 = 500,
\end{align*}
\begin{align*}
    V=\left\lceil K_{1}f^{3/8}\right\rceil, \quad \quad U=\left\lceil \frac{\alpha V}{f^{1/4}}\right\rceil, \quad \quad T=\left\lceil \frac{f^{1/4}}{\alpha}\right\rceil.
\end{align*}
Then $UV\leq 1.001\alpha K_1^2f^{1/2}<f$, as required in Theorem~\ref{thm:1}. By Proposition~\ref{prop2} we have $W(f,T,2)=3+2f^{1/2}/T^2$, and since $UV\geq \alpha K_1^2 f^{1/2}$ and $T^2UV\geq K_1^2 f/\alpha$, it follows that
\begin{align*}
\Delta^{4} = \frac{f^{1/2}W(f,T,2)}{UV} = \frac{3f^{1/2}}{UV}+\frac{2f}{T^{2}UV}\leq \frac{1}{K_1^2}\left(\frac{3}{\alpha}+2\alpha\right) = \frac{2\sqrt{6}}{K_1^2}
\end{align*}
(the choice $\alpha=\sqrt{3/2}$ minimises $3/\alpha+2\alpha$), whence $\Delta\leq 0.0666$.

Write $\varepsilon := E_1+E_2$, and note that $V_0 = UT$ in the notation of Lemma \ref{lem:firsttwo11}; since $UT\leq 1.001K_1f^{3/8}<f$ for $f\geq 10^{22}$, the hypothesis $V_0<f$ of that lemma holds. We have
\begin{align*}
V \leq V_0 \leq (1+\eta)V, \qquad 
\eta := \left( \alpha + \frac{1}{K_1 f^{1/8}} \right) \left( \frac{1}{\alpha} + \frac{1}{f^{1/4}} \right) - 1 \leq 2.5 \times  10^{-5}
\end{align*}
for $f \geq 10^{22}$.
Suppose now, for a contradiction, that $q_2 > 36.88 f^{0.21769}$, and write $q_2 = V_0^{\rho_2}$. Since
\begin{align*}
V_0^{0.5805}\leq \bigl((1+\eta)(K_1f^{3/8}+1)\bigr)^{0.5805}\leq 36.88 f^{0.5805\times 3/8} \leq 36.88 f^{0.21769} < q_2,
\end{align*}
we have $\rho_2 > 0.5805$. On the other hand, by Corollary \ref{lem:firsttwo1} (note that $q_2>q_1\geq 379$),
\begin{align}
\label{eq:rho2bound}
\log\left(\frac{1}{\rho_2}\right)\geq \frac{2(1-\varepsilon)}{3}-\frac{1}{q_1}-\frac{1}{q_2}-\delta_0(q_2)-\frac{1}{2\log^{2} V_0} - \frac{\pi(V_0)}{V}.
\end{align}
We estimate the terms on the right-hand side of \eqref{eq:rho2bound} as follows.
We have $1/q_1\leq 1/379$, and since $q_2>V_0^{0.5805}$,
\begin{align*}
\frac1{q_2}<V_0^{-0.5805}
\quad\text{and}\quad
\delta_0(q_2)\leq\frac{1}{2\log^2\max\{V_0^{0.5805},10^{12}\}};
\end{align*}
moreover, by \cite[Corollary 2]{Ma} we have $q_2\leq 1.821f^{1/4}(\log f)^{3/2}\leq 10^{12}$
for $f\leq 10^{35}$, so that $\delta_0(q_2)=0$ in that range.
By a result of Dusart \cite{Dusart99},
$\pi(x)\leq \frac{x}{\log x}\left(1+\frac{1}{\log x}+\frac{2.51}{\log^2 x}\right)$ for $x\geq 355991$; since $V_0\geq V\geq K_1f^{3/8}\geq 10^{10}$, this gives
\begin{align*}
\frac{\pi(V_0)}{V}\leq (1+\eta)\frac{\pi(V_0)}{V_0}\leq \frac{1+\eta}{\log V_0}\left(1+\frac{1}{\log V_0}+\frac{2.51}{\log^2 V_0}\right).
\end{align*}
All of these bounds are non-increasing in $V_0$, and $V_0\geq K_1f^{3/8}$. Splitting $[10^{22}, 10^{50}]$ into the subranges below, bounding the error terms at the left endpoint of each subrange, and bounding $\varepsilon$ on each subrange by a monotone envelope (the terms of $E_1+E_2$ that increase with $U$ evaluated at the right endpoint, the remaining terms at the left endpoint), we obtain the following table:
\begin{center}
\begin{tabular}{c|c|c|c}
$f$ & $\varepsilon\leq$ & error terms $\leq$ & $\log(1/\rho_2)\geq$ \\\hline
$[10^{22},1.3\times10^{22}]$ & $0.11658$ & $0.04483$ & $0.54411$\\
$[1.3\times10^{22},2\times10^{22}]$ & $0.11660$ & $0.04465$ & $0.54428$\\
$[2\times10^{22},6\times10^{22}]$ & $0.11673$ & $0.04437$ & $0.54447$\\
$[6\times10^{22},3\times10^{23}]$ & $0.11679$ & $0.04366$ & $0.54514$\\
$[3\times10^{23},10^{25}]$ & $0.11716$ & $0.04267$ & $0.54589$\\
$[10^{25},10^{28}]$ & $0.11790$ & $0.04067$ & $0.54739$\\
$[10^{28},10^{33}]$ & $0.11919$ & $0.03726$ & $0.54994$\\
$[10^{33},10^{35}]$ & $0.11900$ & $0.03276$ & $0.55457$\\
$[10^{35},10^{41}]$ & $0.12140$ & $0.03193$ & $0.55380$\\
$[10^{41},10^{50}]$ & $0.12461$ & $0.02823$ & $0.55536$\\
\end{tabular}
\end{center}
(The numerical bounds were verified computationally \cite{code}.)
In every case we conclude that
\begin{align*}
\log\left(\frac{1}{\rho_2}\right)\geq 0.54411 > 0.54387 > \log\left(\frac{1}{0.5805}\right),
\end{align*}
so that $\rho_2 < 0.5805$, a contradiction. Therefore $q_1\leq q_2\leq 36.88f^{0.21769}$,
which completes the proof.
\end{proof}

In Lemma \ref{lem:main2m} below, we will give a variant of Vinogradov's trick which allows an estimation of $m$ similar to Proposition \ref{thm:main2}, but when $q_1$ is bounded from below by a constant. This is based on the following:
\begin{lemma}
\label{lem:large1}
Let $q$ be prime, $\chi$ modulo $q$ be a character of order $3$, $\zeta$ be a primitive third root of unity and $\cS$ be a finite multiset of integers coprime to $q$. Suppose there exists some $i=0,1,2$ such that 
\begin{align*}
|\{ n\in \cS \ : \ \chi(n)=\zeta^{i}\}|\leq K.
\end{align*}
Then 
\begin{align}\label{oban}
\left|\sum_{n\in \cS}\chi(n)\right|\geq \frac{|\cS|-3K}{2}.
\end{align}
\end{lemma}
\begin{proof}
Since the bound in \eqref{oban} is decreasing with $K$ we may assume 
\begin{equation*}
|\{ n\in \cS \ : \ \chi(n)=\zeta^{i}\}|=K.
\end{equation*}
Let 
$S=\sum_{n\in \cS}\chi(n)$,
and decompose 
\begin{align*}
S=\sum_{j=0}^{2}|\cA_j|\zeta^j,
\end{align*}
where $\cA_j=\{ n\in \cS \ : \ \chi(n)=\zeta^{j}\}$.
We have 
\begin{align*}
|S|=|\zeta^{-i}S|\geq \left|\sum_{j\neq i}|\cA_j|\zeta^{j-i}\right|-K.
\end{align*}
Hence there exists integers $N_1,N_2$ satisfying 
$N_1+N_2=|\cS|-K$,
such that 
\begin{align*}
|S|\geq \left|N_1e^{2\pi i/3}+N_2e^{-2\pi i /3} \right|-K.
\end{align*}
The result follows since 
$$\left|N_1e^{2\pi i/3}+N_2e^{-2\pi i /3} \right|\geq (N_1+N_2)\cos(\pi/3) = \frac{|\cS|-K}{2},$$
so that
$$|S| \geq \frac{|\cS|-K}{2}-K = \frac{|\cS|-3K}{2}.$$
\end{proof}
\begin{lemma}
\label{lem:main2m}
Assume $10^{22} \leq f \leq 10^{50}$, let $\omega$ be a cube root of unity, and suppose
$q_1 \ge 379$.
Then there exists a positive integer $m$ satisfying
\begin{align} \label{eq: prop23-cond-m}
\chi(m)=\omega, \quad (m,q_1 q_2)=1, \quad m \leq 14 f^{1/3}(\log{f})^{1/2}.
\end{align}
\end{lemma}
\begin{proof}
For $T,U,V \in \N$, define the sets
\begin{align*}
S := & \{(t,u,v) \in \N^3 \mid \ t \leq T , \ u \leq U, v \leq V \}, \\
S^* := & \{(t,u,v) \in S \mid (v-ut,q_1q_2) = 1\},
\end{align*}
and the sums 
\begin{align*}
\Sigma_S:=\sum_{(t,u,v) \in S}\chi(v-ut) \quad \text{and} \quad
\Sigma_{S^{*}}:=\sum_{(t,u,v) \in S^*}\chi(v-ut).
\end{align*}
By the triangle inequality 
\begin{equation} \label{eq: lem23-triangle-ineq}
|\Sigma_S-\Sigma_{S^{*}}|\leq \sum_{(t,u,v) \in S \setminus S^*} 1.
\end{equation}
Since $q_1,q_2$ are prime
\begin{equation*} 
\sum_{(t,u,v) \in S \setminus S^*} 1 \leq \sum_{\substack{(t,u,v) \in S \\ q_1 \mid v-ut}}1+\sum_{\substack{(t,u,v) \in S \\ q_2 \mid v-ut}}1 \leq TUV \left(\frac{1}{q_1}+\frac{1}{q_2} + \frac{2}{V}\right) =: \delta\, TUV.
\end{equation*}
Let us set:
\begin{align*}
r=3, \quad T=\left\lceil f^{1/6}\right\rceil, \quad U=\left\lceil \frac{V}{T}\right\rceil, \quad \theta=0.53, \quad \text{and } V=\lfloor 13 f^{1/3} \log^{1/2} f \rfloor,
\end{align*}
then, since $q_2>q_1\ge 379$ and $f \geq 10^{22}$, we get $q_2 \geq 383$ and thus $\delta \leq 5.3 \times 10^{-3}$.

Invoking Theorem \ref{thm:1} with the parameters chosen as above, we obtain
\begin{align*}
|\Sigma_{S}|\leq TUV(E_1+E_2)<0.49 \times  TUV,
\end{align*}
where $E_1+E_2$ is evaluated exactly \cite{code}, with all floors and ceilings retained; it is decreasing on $[10^{22},10^{50}]$, with maximum $0.48296$ at $f=10^{22}$. Thus by \eqref{eq: lem23-triangle-ineq},
\begin{align*}
|\Sigma_{S^{*}}|< (0.49 + \delta) TUV.
\end{align*}
Consider the multiset $\mathbf{S^*} := \{v-ut \mid (t,u,v) \in S^*\}$. Every element of $\mathbf{S^*}$ is non-zero, since $(0,q_1q_2)\neq1$, and satisfies $|v-ut|\leq\max\{V,UT\}<f$, hence is coprime to $f$. Assume, for a contradiction, that $|\{m \in \mathbf{S^*} \mid \chi(m) = \omega \}| = 0$. By Lemma \ref{lem:large1}, we have
$$ \left| \Sigma_{S^*} \right| = \left| \sum_{m \in \mathbf{S^*}} \chi(m) \right| \geq \frac{|\mathbf{S^*}|}{2}, \quad \text{implying} \quad  |\mathbf{S^*}|=|S^*| \leq 2 (0.49 + \delta) \, TUV. $$
On the other hand,
$$ |S^*| = |S|-|S\setminus S^*| \geq (1-\delta) TUV,$$
from where we get
\begin{align*}
1- \delta \leq 2 (0.49 + \delta) \quad \text{implying} \quad \frac1{150} \leq \delta,
\end{align*}
which leads to a contradiction. Thus, for any cube root of unity $\omega$, there is a non-zero integer $m$, $|m| \leq \max\{TU, V\}$, coprime to $q_1 q_2$ and satisfying $\chi(m)=\omega$. Since $\chi$ has order $3$, we have $\chi(-1)=1$, so replacing $m$ by $|m|$ preserves the value of $\chi(m)$. Finally, since $V \leq UT \leq 14 f^{1/3} \log^{1/2} f$ for $f \geq 10^{22}$, $|m|$ is a positive integer satisfying \eqref{eq: prop23-cond-m}.
\end{proof}

\section{Proof of Theorem \ref{T:2}}
\label{sec:II}
Suppose $f\geq 10^{22}$.
By Proposition \ref{fixedprop14}, we may assume $q_1\geq 379$ and $q_2 \geq 383$.
Recall that we aim at showing that, with $\omega:=\chi^{-1}(q_2)$,
there exists a positive integer $m$ such that $(m,q_1q_2)=1$, and $\chi(m)=\omega$, and
\begin{align}
\label{eq:cond2}
3q_1 q_2 m<f,
\end{align}
\begin{align}
\label{eq:cond3}
10 q_1^2 q_2<f.
\end{align}
By \cite{Lezowski} we may further assume $f\leq 10^{50}$, so that
Theorem \ref{thm:main1} applies and yields $q_1\leq q_2\leq 36.88f^{0.21769}$. In particular,
\begin{align*}
10q_1^2q_2\leq 10\cdot 36.88^3 f^{3\times 0.21769}<f,
\end{align*}
since $f^{0.34693}>501670$, which settles \eqref{eq:cond3}. We split the proof of \eqref{eq:cond2} into two cases.
\vspace{1ex}

\noindent
{\bf Case I:}  $q_1 \leq 1.82 f^{1/8}$.
By Lemma \ref{lem:main2m},
$$m\leq 14 f^{1/3}\log^{1/2}{f},$$
and hence
\begin{align*}
3q_1q_2m\leq 3\cdot 1.82f^{1/8}\cdot 36.88f^{0.21769}\cdot 14f^{1/3}\log^{1/2}f
< 2820\, f^{0.677}\log^{1/2}f < f,
\end{align*}
since $2820\log^{1/2}f\leq f^{0.323}$ for $f\geq 10^{22}$.
\vspace{1ex}

\noindent
{\bf Case II:} $q_1 > 1.82 f^{1/8}$.
By Proposition \ref{thm:main2} we have
$$m\leq 86 f^{5/16},$$
and hence
\begin{align*}
3q_1 q_2 m \leq 3\cdot 36.88^2 f^{2\times 0.21769}\cdot 86 f^{5/16} < 350915\, f^{0.74789} < f,
\end{align*}
since $f^{0.25211} > 350915$ for $f\geq 10^{22}$. This completes the proof of the theorem.

\section{Numerical verification}
\label{numerics}
In this section we describe an algorithm to enumerate primes $f$ with no small cubic non-residues modulo $f$, and its application to prove
the following result, completing the proof of Theorem~\ref{T:big}.
See \cite{code} for the source code.

\begin{theorem}\label{T:compute}
There are no cyclic cubic norm-Euclidean number fields with conductor
$f\in(2\times10^{14},10^{22})$.
\end{theorem}

The algorithm is similar to methods of enumerating pseudosquares and
pseudocubes (see \cite{Sorenson} and the references therein), the main
differences being that our sieve region is an ellipse rather than an
interval, and that the sieve conditions depend on both coordinates of a
point $(x,y)$. Over most of the range we employ Bernstein's
``doubly-focused enumeration'' strategy \cite{Bernstein}, as described
below.

We begin with some basic notation and lemmas.
Let $\omega=\frac{-1+\sqrt{-3}}{2}$. For $\alpha,\pi\in\Z[\omega]$
with $\pi\nmid3$ prime, recall the cubic residue symbol
$\left(\frac{\alpha}{\pi}\right)_3\in\bigl\{0,1,\omega,\omega^{-1}\bigr\}$
defined by
\begin{equation*}
\left(\frac{\alpha}{\pi}\right)_3
\equiv\alpha^{\frac{\pi\overline{\pi}-1}{3}}
\pmod\pi.
\end{equation*}
\begin{lemma}\label{lem:2and3}
Let $p\equiv1\pmod3$ be a prime such that $2$ and $3$ are cubic residues
modulo $p$. Then $p$ can be written uniquely in the form $x^2+243y^2$,
where $x,y\in\Z$, $x\equiv 3y-1\pmod6$, and $y>0$.
\end{lemma}
\begin{proof}
By \cite[Prop.~9.6.2]{IR}, $2$ is a cubic residue modulo $p$ if and only
if $p$ is of the form $C^2+27D^2$ for $C,D\in\Z$, and this representation
is unique up to signs. We choose the signs such that $C\equiv-1\pmod3$ and
$D>0$, and write $\pi=C+3D\sqrt{-3}$, which is a prime of $\Z[\omega]$
satisfying $\pi\overline{\pi}=p$. By \cite[Ch.~9, Thm.~$1'$]{IR}, $3$
is a cubic residue modulo $p$ if and only if
\begin{align*}
1&=\left(\frac{3}{\pi}\right)_3=\left(\frac{-\omega^2(1-\omega)^2}{\pi}\right)_3
=\omega^{\frac23(p-1)}\omega^{\frac43(C+3D+1)}\\
&=\omega^{\frac23(C^2+2C+1+6D+27D^2)}=\omega^{\frac23(C+1)^2}\omega^{4D+18D^2}
=\omega^D.
\end{align*}
Hence $3\mid D$, and writing $C=x$, $D=3y$, we have $p=x^2+243y^2$, uniquely.
Finally note that since $p$ is odd, $x$ and $y$ must have opposite parity,
and hence $x\equiv3y-1\pmod6$.
\end{proof}

\begin{lemma}\label{lem:q}
Let $p=x^2+243y^2=\pi\overline{\pi}$ as in Lemma~\ref{lem:2and3},
and let $q>3$ be a prime number.
\begin{enumerate}
\item[(i)]
If $q\equiv1\pmod3$ then, as elements of $\F_q$,
\begin{equation*}
\bigl(p(x-9ys)\bigr)^{\frac{q-1}{3}}=
\begin{cases}
0&\text{if }\left(\frac{q}{\pi}\right)_3=0,\\
1&\text{if }\left(\frac{q}{\pi}\right)_3=1,\\
\frac{-1+s}{2}&\text{if }\left(\frac{q}{\pi}\right)_3=\omega,\\
\frac{-1-s}{2}&\text{if }\left(\frac{q}{\pi}\right)_3=\omega^{-1},
\end{cases}
\end{equation*}
where $s\in\F_q$ is any square root of $-3$.
\item[(ii)]
If $q\equiv2\pmod3$ then
$\bigl(x+9y\sqrt{-3}\bigr)^{\frac{q^2-1}{3}}=\left(\frac{q}{\pi}\right)_3$
in $\F_q\bigl[\sqrt{-3}\bigr]$.
\end{enumerate}
In particular, $q$ is a cubic residue modulo $p$ if and only if
$\bigl(p(x-9ys)\bigr)^{\frac{q-1}{3}}=1$ in case (i),
and $\bigl(x+9y\sqrt{-3}\bigr)^{\frac{q^2-1}{3}}=1$ in case (ii).
\end{lemma}
\begin{proof}
Suppose $q\equiv1\pmod3$, and
write $q=\Pi\overline{\Pi}$, where $\Pi\in\Z[\omega]$
with $\Pi\equiv2\pmod{3}$. Then by cubic reciprocity \cite[Ch.~9, Thm.~1]{IR},
\begin{equation*}
\left(\frac{q}{\pi}\right)_3
=\left(\frac{\Pi}{\pi}\right)_3
\left(\frac{\overline{\Pi}}{\pi}\right)_3
=\left(\frac{\pi}{\Pi}\right)_3
\left(\frac{\pi}{\overline{\Pi}}\right)_3.
\end{equation*}
Applying \cite[Prop.~9.3.4]{IR}, this becomes
$\left(\frac{\pi\overline{\pi}^2}{\Pi}\right)_3
=\left(\frac{p\overline{\pi}}{\Pi}\right)_3$. Fixing a choice of $s\in\Z$
with $s^2\equiv-3\pmod{q}$, we may swap $\Pi$ and $\overline{\Pi}$
if necessary to assume that $\sqrt{-3}\equiv s\pmod{\Pi}$. Hence we
have $\overline{\pi}=x-9y\sqrt{-3}\equiv x-9ys\pmod\Pi$, and the claim
follows from the definition of the cubic residue symbol.

Suppose now that $q\equiv2\pmod3$. Then by cubic reciprocity we have
$\left(\frac{q}{\pi}\right)_3=\left(\frac{\pi}{q}\right)_3$, and again by the
definition of the symbol this has image
$\bigl(x+9y\sqrt{-3}\bigr)^{\frac{q^2-1}{3}}$ in
$\F_q\bigl[\sqrt{-3}\bigr]$.
\end{proof}

Given a range $(B_1,B_2)$, we can
use these lemmas to quickly enumerate candidate primes $f\in(B_1,B_2)$
with no small cubic non-residues mod $f$, as follows.
Using Lemma~\ref{lem:q} we precompute the values of
$\left(\frac{q}{\pi}\right)_3$, where $\pi=x+9y\sqrt{-3}$,
for every prime $q\in(3,541]$ and all $(x,y)\in\{0,\ldots,q-1\}^2$.
We then pick out the pairs for which the symbol is $1$ for all small $q$,
and use the Chinese remainder theorem to combine them, as follows.
Let $M=M_1M_2$ be a product of the first several primes with
$(M_1,M_2)=1$, and for a given $y$, call $x\in\Z/M\Z$
\textit{admissible} if $(x,y)$ reduces to an admissible pair modulo every
prime dividing $M$. The admissible $x$ are precisely the values
\begin{equation*}
x\equiv x_1'+x_2'\pmod{M},
\qquad
x_i'=\bigl(x_i\overline{M}_{3-i}\bmod M_i\bigr)M_{3-i},
\end{equation*}
where $\overline{M}_i$ denotes a multiplicative inverse of
$M_i\pmod{M_{3-i}}$, and $x_i$ runs over the residues for which
$(x_i,y)\in(\Z/M_i\Z)^2$ is admissible. Since
$0\leq x_1'+x_2'<2M$, every admissible
$x\in[-X,X]$, where $X=\sqrt{B_2-243y^2}$, has the form
\begin{equation}\label{eq:dfe}
x = x_1'+x_2'-jM \quad\text{for some integer } j.
\end{equation}
Note also that any prime $f$
omitted by this sieve has least non-residue $q_1\mid M$, so
the omission is justified by Theorem~\ref{T:1} provided that
$B_1\geq f_0(q)$ for the largest prime $q\mid M$.

For $(B_1,B_2)=\bigl(2\times10^{14},3.26\times10^{16}\bigr)$ we choose
$M_1=2\cdot3\cdot5$ and $M_2=7$
(so that $B_1\geq f_0(7)=10^{14}$), and precompute the lists
$\{x_i'\}$ for every residue $y\bmod{M_i}$. In this range
$M=210$ is much smaller than $X$ for most $y$, and for each pair
$(x_1',x_2')$ we simply traverse the class
$x\equiv x_1'+x_2'\pmod{M}$ in $[-X,X]$ in steps of $M$.

For $(B_1,B_2)=\bigl(3.26\times10^{16},10^{22}\bigr)$ we take
$M$ to be the product of the primes up to $41$, split as
\begin{equation*}
M_1=2\cdot3\cdot11\cdot13\cdot23\cdot31\cdot37,
\qquad
M_2=5\cdot7\cdot17\cdot19\cdot29\cdot41
\end{equation*}
(so that $B_1=f_0(41)$). In this range $M\approx3\times10^{14}$
exceeds $2X$, so each admissible residue class contains at most one
point of $[-X,X]$, and \eqref{eq:dfe} can only occur with
$j\in\{0,1,2\}$; on the other hand, with $M_i\approx10^7$ it is no
longer feasible to precompute the lists $\{x_i'\}$ for every residue
$y\bmod{M_i}$. Instead, we factor
$M_1=(2\cdot3\cdot31\cdot37)(11\cdot13\cdot23)$ and
$M_2=(5\cdot19\cdot41)(7\cdot17\cdot29)$, and generate the two
lists afresh for each $y$ by applying the identity above to these
factorizations, so that precomputed lists are needed only for the four
smaller factors. We then sort $\{x_1'\}$ and $\{x_2'\}$ and extract
all sums $x_1'+x_2'$ lying in
$[0,X]\cup[M-X,M+X]\cup[2M-X,2M)$ --- that is, all solutions of
\eqref{eq:dfe} --- by three linear merges with sliding pointers, at a
total cost proportional to the lengths of the lists plus the number of
solutions. This is Bernstein's doubly-focused enumeration, transposed
from an interval to our elliptical region.

We process $f=x^2+243y^2$ for a fixed pair $(x,y)\in\Z^2$ as follows:
\begin{enumerate}
\item We search for a cubic non-residue $q_1\le199$, and if successful we reject $f$ if it exceeds the corresponding threshold in Table~\ref{tab:q1thresh}.
\item If no $q_1\le199$ is found, we test to see if $\gcd(x,y)>1$ or $f$ is a cube, and reject $f$ if so; otherwise we continue the search for $q_1$ up to $541$.
\item If $q_1\le541$ then we search for prime cubic non-residues $q_2,m$ such that $q_1<q_2<m\le 541$ and $q_2m$ is a cubic residue, and apply Criterion~\ref{criterion: expl-non-Eucl}.
\item Any candidates that remain at this point (either because we could not find suitable $q_1,q_2,m\le541$ or the criterion was not met) are printed.
\end{enumerate}
We also reject any candidate $f$ that is revealed to be composite by trial division along the way. However, we never apply a general primality test, since it is quicker to complete the steps above, even when $f$ is \emph{de facto} composite.

We ran this algorithm on a desktop PC with 6 cores (Intel Core i7-8700). The running time was approximately 15 hours for $(B_1,B_2) = (2\times10^{14},3.26\times10^{16})$, and
13 days for $(B_1,B_2) = (3.26\times10^{16},10^{22})$.
Neither run printed any exceptional values of $f$. This completes the proof of Theorem~\ref{T:compute}.

Finally, we have:

\begin{proof}[Proof of Theorem~\ref{T:big}]
In light of~\cite{Lezowski} we may assume $2\times 10^{14}\leq f\leq 10^{50}$.
Combining Theorems~\ref{T:1},~\ref{T:2}, and~\ref{T:compute}
gives the result.
\end{proof}

\section*{Acknowledgements}
We thank Olivier Ramar\'{e} for discussions on this topic.
We also wish to thank the UNSW Canberra Rector's Visiting Fellowship, which enabled KM to visit BK and TT in February and March 2019. 

\section*{Conflict of Interest Statement}
On behalf of all authors, the corresponding author states that there is no conflict of interest.

\end{document}